\documentclass[11pt]{amsart}
\usepackage{latexsym, amsmath,amssymb}
\usepackage{hyperref}
\usepackage{tikz}
\usetikzlibrary{calc,trees,shadows,positioning,arrows,chains,shapes.geometric,
decorations.pathreplacing,decorations.pathmorphing,shapes,
matrix,shapes.symbols,patterns,intersections,fit}

\theoremstyle{plain}
\newtheorem{thm}{Theorem}

\newtheorem{lem}{Lemma}
\newtheorem{prop}{Proposition}
\newtheorem{cor}{Corollary}

	\newcommand{\D}{\mathcal{D}}
	\newcommand{\eps}{\varepsilon}
		\newcommand{\ls}{\lesssim}
				\newcommand{\gs}{\gtrsim}
			\newcommand{\la}{\lambda}
			\newcommand{\1}{{\rm 1\hspace*{-0.4ex}%
					\rule{0.1ex}{1.52ex}\hspace*{0.2ex}}}
\begin{document}

\title[Strichartz estimates for orthonormal systems]{Strichartz estimates for orthonormal systems on compact manifolds} 

\begin{abstract}
We establish new Strichartz estimates for orthonormal systems on compact Riemannian manifolds in the case of wave, Klein-Gordon and fractional Schr\"odinger equations. Our results generalize the classical (single-function) Strichartz
estimates on compact manifolds by Kapitanski \cite{kap89}, Burq-G\'erard-Tzvetkov \cite{MR2058384}, Dinh \cite{dinh17}, and extend the Euclidean orthonormal version by Frank-Lewin-Lieb-Seiringer \cite{MR3254332}, Frank-Sabin \cite{MR3730931}, Bez-Lee-Nakamura  \cite{BLN21}. On the flat torus, our new results for the Schr\"odinger equation cover  prior work of Nakamura \cite{MR4068269}, which exploits the dispersive estimate of Kenig-Ponce-Vega \cite{MR1101221}. We achieve sharp results on compact manifolds by  combining the frequency localized  dispersive estimates for small time intervals with the duality principle due to Frank-Sabin.  We construct examples to show these results can be saturated on the sphere, and we can  improve them on the flat torus by using Bourgain-Demeter's decoupling theorem to obtain new   decoupling inequalities for certain non-smooth hypersurfaces. As an application, we obtain the well-posedness of infinite systems of dispersive equations with Hartree-type nonlinearity.
 \end{abstract}
	



\author{Xing Wang, An Zhang and Cheng Zhang}

\address{School of Mathematics, Hunan University, Changsha, HN 410012, China}
\email{xingwang@hnu.edu.cn}

\address{School of Mathematical Sciences, Beihang University, Beijing, BJ 100191, China} 
\email{anzhang@buaa.edu.cn} 

	\address{Yau Mathematical Sciences Center,
	Tsinghua University,
	Beijing, BJ 100084, China}
\email{czhang98@tsinghua.edu.cn}



\maketitle

\section{Introduction}	\label{int}
Let $d\ge 1$.   Let $(M,g)$ be a $d$-dimensional smooth compact Riemannian manifold without boundary. Let $\Delta_g$ be the Laplacian-Beltrami operator on $M$. Let $\Delta=-\Delta_g$. Let $e^{itP}f$ denote the solution to the initial value problem
\begin{equation*}
\left\{\begin{array}{l}
i\partial_t u+Pu=0, \quad (x,t)\in M\times \mathbb R,\\
u(\cdot,0)=u_0\,.
\end{array}\right.
\end{equation*} 
A wide class of dispersive equations have this form, such as the Schr\"odinger equation $P=\Delta$, the fractional Schr\"odinger equation $P=\Delta^{\alpha/2}$ ($\alpha\ne0,1$), the wave equation $P=\sqrt\Delta$, and the Klein-Gordon equation $P=\sqrt{1+\Delta}$. 

In quantum mechanics, a system of $N$ independent fermions is described by a collection of $N$ orthonormal functions $f_1,...,f_N$ in $L^2$. So functional inequalities that incorporate a significant number of orthonormal functions are highly valuable for the mathematical analysis of large-scale quantum systems. The inequalities have applications to the Hartree equation modeling infinitely many fermions in a quantum system, see Chen–Hong–Pavlovic \cite{CHP17,CHP18}, Frank–Sabin \cite{MR3730931}, Lewin–Sabin \cite{MR3304272,MR3270166} and Sabin \cite{sabin14}. The idea in this line of investigation is to generalize the classical inequalities for a single-function input to an orthonormal system. In the pioneering work of Lieb-Thirring \cite{LT}, they first established such an extension of the Gagliardo–Nirenberg–Sobolev inequality. In the recent work of Frank-Lewin-Lieb-Seiringer \cite{MR3254332} for the Schr\"odinger propagator $e^{it\Delta}$, they proved a generalization of the well known Strichartz estimates for systems of orthonormal functions in $L^2(\mathbb{R}^d)$. Later, Frank-Sabin \cite{MR3730931}, Bez-Hong-Lee-Nakamura-Sawano \cite{MR3985036}, Bez-Lee-Nakamura  \cite{BLN21}, Feng-Mondal-Song-Wu \cite{FMSW24} investigated a wide class of dispersive equations and established the Strichartz estimates for systems of orthonormal functions on the Euclidean space.   To our best knowledge, there are only a few results concerning such generalizations on compact Riemannian manfiolds. Frank-Sabin \cite{FS} established the spectral cluster bounds for orthonormal systems on compact manifolds, and recently Ren-Zhang \cite{renzhang24} obtained some  improvements on non-positively curved manifolds. Nakamura \cite{MR4068269} studied the Strichartz estimates on the flat torus with orthonormal system input and obtained sharp estimates in certain sense. In this paper, we provide substantial progress in this direction, extending the  orthonormal Strichartz estimates to general compact Riemannian manifolds.

The classical (single-function) Strichartz estimates in the Euclidean space date back to the seminal paper of Strichartz \cite{MR512086}.  See also  Ginibre-Velo \cite{MR1151250}, Keel-Tao \cite{MR1646048}  and references therein. In the case of compact manifolds, Kapitanski \cite{kap89}, Burq-G\'erard-Tzvetkov \cite{MR2058384},  Dinh \cite{dinh17} obtained Strichartz estimates for the wave, the Schr\"odinger and the fractional Schr\"odinger equations respectively. See also Cacciafesta-Danesi-Meng \cite{cdm23} for the Dirac equations. In the case of the torus, see e.g.  the celebrated work of Bourgain-Demeter \cite{MR3374964}. In this paper, we shall investigate  their generalizations to orthonormal  systems.

In the following, let $I\subset\mathbb{R}$ be a fixed bounded interval of length $|I|\approx 1$. We fix $(e_k)_k$ to be an orthonormal eigenbasis in $L^2(M)$ associated with the eigenvalues $(\la_k)_k$ of $\sqrt{\Delta}$. Here $0=\la_0<\la_1\le \la_2\le ...$ are arranged in increasing order and we account for multiplicity. We define the Fourier coefficient  $\hat f(k)=\langle f,e_k\rangle$ for  $f\in L^2(M)$. We call $(f_j)_j$  an  orthonormal system in a Hilbert space $\mathcal{H}$ if the functions $f_j$ are orthonormal in $\mathcal{H}$. We investigate the Strichartz estimates of the form 
\begin{equation}\label{eq0}
	\Big\|\sum_j  \nu_j |e^{itP} f_j|^2\Big\|_{L_t^{p/2}L_x^{q/2}(I\times M)} \lesssim N^{\sigma}
	\|\nu\|_{l^\beta}\end{equation}
for all orthonormal systems $(f_j)_j$ in $L^2(M)$ with $supp\, \hat f_j\subset \{k:\la_k\le N\}$, and all sequences $\nu=(\nu_j)_j\in l^\beta$. The estimates are independent of the choice of the eigenbasis $(e_k)_k$. The main goal in the line of investigation is to determine the optimal range of $\beta$ for a fixed exponent $\sigma$. A natural choice of $\sigma$ is just the one in the classical (single-function) case. In this case, \eqref{eq0} trivially holds with $\beta=1$ by Minkowski inequality, while the question is to determine the largest $\beta$ by exploiting the orthogonality between the functions. On the other hand, it is also interesting to determine  the optimal range of $\beta$ for any fixed $\sigma$. We shall establish sharp estimates in the form of \eqref{eq0} by combining the frequency localized  dispersive estimates for small time intervals with the duality principle due to Frank-Sabin. In  comparison, Nakamura \cite{MR4068269} and Bez-Lee-Nakamura  \cite{BLN21} exploit stronger frequency global  dispersive estimates in the flat case, as in the prior work of Kenig-Ponce-Vega \cite{MR1101221}.

$ $

\noindent\textbf{Notations.} Throughout this paper, $X\ls Y$ means $X\le CY$  for some positive constants $C$. If the constant depends on $\eps>0$, we denote $X\ls_\eps Y$. If $X\ls Y$ and $Y\ls X$, we denote $X\approx Y$.

	Now, we introduce our main  results on the fractional Schr\"odinger, the wave and the Klein-Gordon equations on compact manifolds.
\subsection{Fractional Schr\"odinger equations}
\begin{figure}[h]
	\centering
	\includegraphics[width=0.7\textwidth]{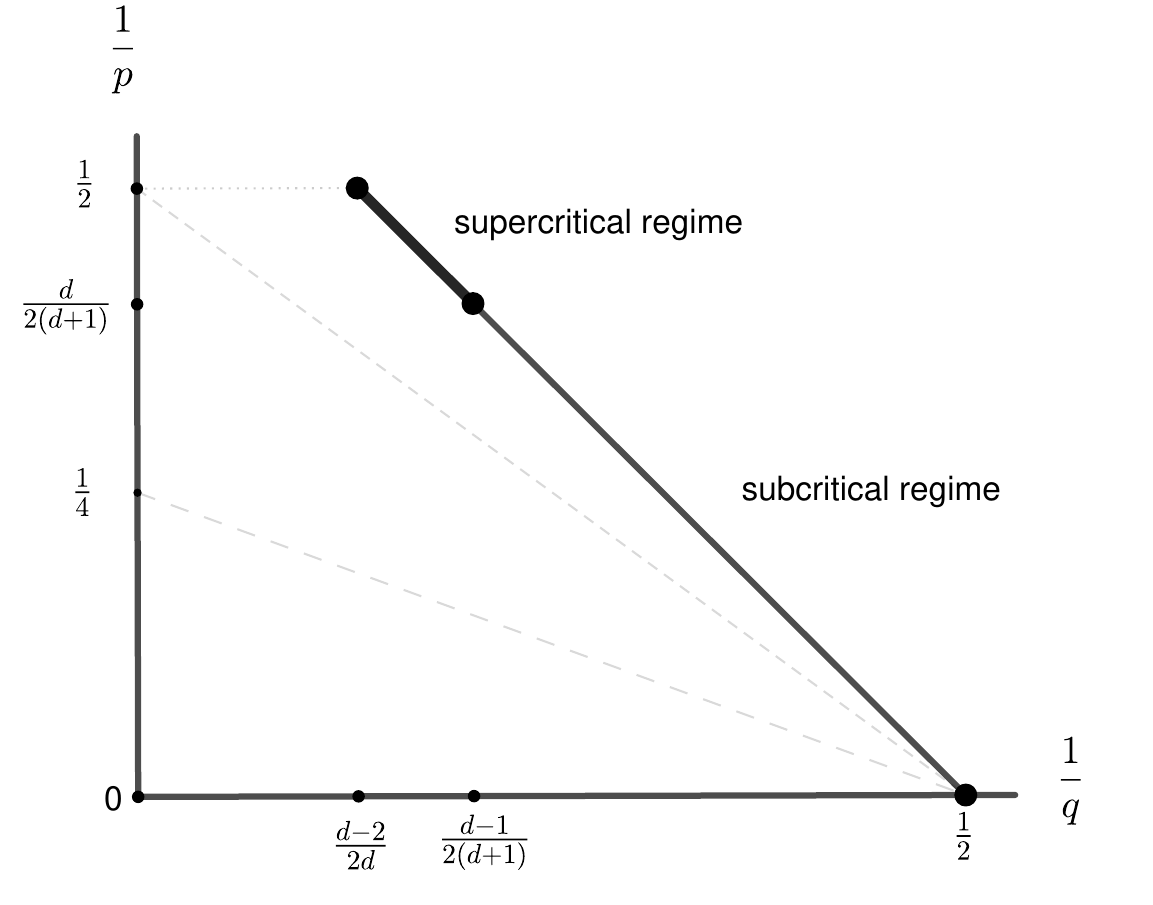}
	\caption{Sharp Schr\"odinger admissible pairs}
	\label{fig1}
\end{figure}

	 Suppose $p\ge 2,\ q<\infty$ and $\frac1p=\frac d2(\frac12-\frac1q)$.  We divide these sharp Schr\"odinger admissible pairs $(p,q)$  into four groups. See Figure \ref{fig1}.
	
	(i) Subcritical regime: $d\ge1$, $2\le q<\frac{2(d+1)}{d-1}$
	
	(ii) Critical point: $d\ge2$, $q=\frac{2(d+1)}{d-1}$
		
	(iii) Supercritical regime: $d=2$, $\frac{2(d+1)}{d-1}< q<\infty$ or $d\ge3$, $\frac{2(d+1)}{d-1}< q< \frac{2d}{d-2}$

	(iv) Keel-Tao endpoint: $d\ge3$, $q=\frac{2d}{d-2}$.

\begin{thm}\label{thm1}
Let $d\ge 1,\,\alpha \in(0,\infty)\setminus\{1\},\ N\ge10$. Suppose $p\ge 2,\ q<\infty$ and $\frac1p=\frac d2(\frac12-\frac1q)$. Let
\begin{align}\label{gam0}
	\sigma_0(\alpha)=\begin{cases}
		2/p,\ \ \ \ \ \ \ \ \ \ \ \ \alpha>1\\
		2(2-\alpha)/p,\ \ \ \alpha\in(0,1).
	\end{cases}
\end{align}
Then 
\begin{equation}\label{str}
	\Big\|\sum_j  \nu_j |e^{it\Delta^{\alpha/2}} f_j|^2\Big\|_{L_t^{p/2}L_x^{q/2}(I\times M)} \lesssim 
	N^{\sigma_0}\|\nu\|_{l^\beta}\end{equation}
holds for all orthonormal systems $(f_j)_j$ in $L^2(M)$ with $supp\, \hat f_j\subset \{k:\la_k\le N\}$, and all sequences $\nu=(\nu_j)_j\in l^\beta$, and the following $\beta$ with respect to the pairs $(p,q)$  in the four groups:

{\rm(i)} Subcritical regime: $\beta\le \frac{d}{d-2/p}=\frac{2q}{q+2}$

{\rm(ii)} Critical point: $\beta<p/2$

{\rm(iii)} Supercritical regime:  $\beta< p/2$

{\rm(iv)} Keel-Tao endpoint: $\beta=1$.
\end{thm}

\begin{figure}[h]
	\centering
	\includegraphics[width=0.8\textwidth]{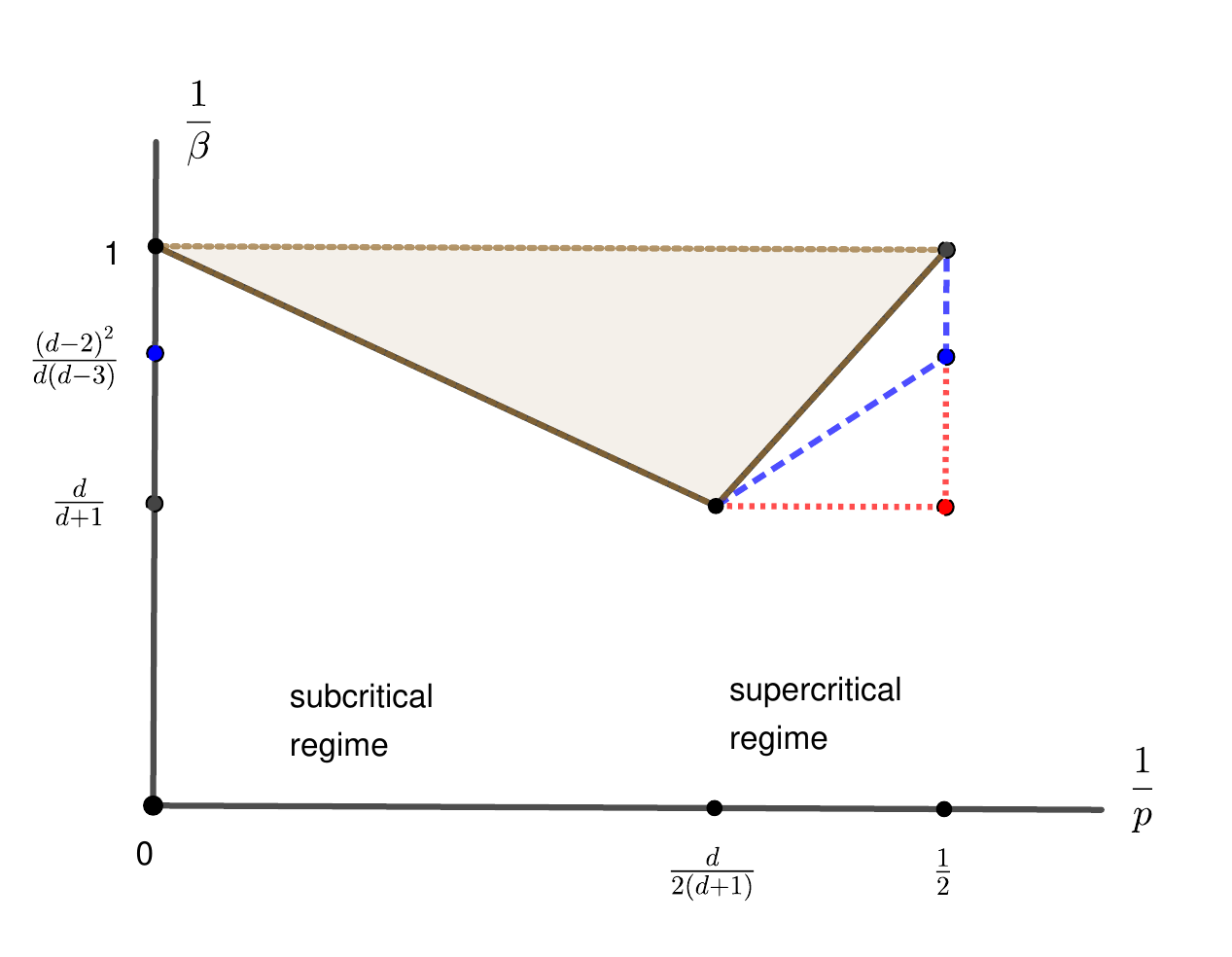}
	\caption{Fractional Schr\"odinger equations with $\alpha\ge2$. General manifolds (shaded triangle, Theorem \ref{thm1}). Improvement on the flat torus  based on Bourgain-Demeter's decoupling theorem (blue, Corollary \ref{torussup}). Conjecture on the flat torus based on the Discrete Restriction Conjecture (red).}
	\label{fig3}
\end{figure}

When $\alpha>1$, the ranges of $\beta$  are  sharp   in the subcritical regime on any compact manifold by the necessary condition \eqref{nec},  and sharp in the supercritical regime on the sphere by the necessary condition \eqref{nec22}.  As we will observe in the forthcoming Corollary \ref{torussup}, the range of $\beta$ in the supercritical regime  can be improved on the flat torus. See Figure \ref{fig3}. This suggests that the optimal range of $\beta$ in the supercritical regime should be sensitive to the geometry of the manifold. It is interesting to determine the optimal range of $\beta$ in the supercritical regime on the manifolds under certain geometric assumptions, such as the hyperbolic manifolds. Furthermore, when $0<\alpha<1$, we shall see from the proof that the estimate \eqref{str} indeed holds on longer intervals of length  $|I|\approx N^{1-\alpha}$, and then the ranges of $\beta$  are sharp in the subcritical and supercritical regimes on any compact manifold  by the necessary conditions \eqref{nec} and \eqref{nec221}.

The exponent $\sigma_0/2$ is exactly the Sobolev exponent in the  classical (single-function) Strichartz estimates by Burq-G\'erard-Tzvetkov \cite{MR2058384} and Dinh \cite{dinh17}. The optimality of $\beta$ only makes sense when the exponent of $N$ in \eqref{str} is fixed, so we fix $\sigma_0/2$ to be the one in the single-function case. Moreover, Nakamura \cite[Theorem 1.5]{MR4068269} obtained the same  estimates in the subcritical regime for the Schr\"odinger propagator $e^{it\Delta}$ on the flat torus. Furthermore, we may expect to raise the exponent of $N$ to increase the range of $\beta$.  See Theorem \ref{thm2}.
	
Next, we obtain  new (single-function) Strichartz estimates on the flat torus for fractional Schr\"odinger equations. We achieve this by establishing a new decoupling inequality for the hypersurface $(\xi,|\xi|^\alpha) \in \mathbb{R}^{d+1}$. See  Theorem \ref{dec0}. We shall exploit Bourgain-Demeter's $l^2$ decoupling theorem for the paraboloid \cite[Theorem 1.1]{MR3374964}. On the flat torus $\mathbb{T}^d=\mathbb{R}^d/\mathbb{Z}^d$, it is standard to use the Fourier coefficients  $\hat f(k)=\langle f,e^{2\pi ik\cdot x}\rangle$ associated to the orthonormal basis $\{e^{2\pi ik\cdot x}\}_{k\in \mathbb{Z}^d}$. For simplicity, we  only discuss $\alpha\ge2$ here. See the forthcoming Theorem \ref{fracstr} and Theorem \ref{fracstr2} for a complete discussion on $\alpha \in(0,\infty)\setminus\{1\}$.
\begin{thm}\label{fracstr0}
	Let $d\ge 1,\,\alpha\ge2,\ N\ge10$.  Let $f\in L^2(\mathbb{T}^d)$ with $supp \hat f\subset [-N,N]^d$.  Suppose $p\ge 2,\ q<\infty$ and $\frac1p=\frac d2(\frac12-\frac1q)$.  Then  for all $\eps>0$, 
	\begin{equation}\label{dec01}
		\|e^{it\Delta^{\alpha/2}}f\|_{L^p_{t}L^q_x(\mathbb{T}^{d+1})}\ls_\eps N^{\sigma_1}\|f\|_{L^2(\mathbb{T}^d)}
	\end{equation}
	where 
	\begin{equation}\label{sigma01}
		\sigma_1(\alpha)=
			\max\Big\{0,\frac d2-\frac{d+2}q\Big\}+\eps.
	\end{equation}
\end{thm}
The case $\alpha=2$ is due to Bourgain-Demeter \cite[Theorem 2.4]{MR3374964}, and the case $\alpha>2$ is new. In particular, at the Keel-Tao endpoint $(p,q)=(2,\frac{2d}{d-2})$ we have 
\begin{equation}\label{keeltao}	\|e^{it\Delta^{\alpha/2}}f\|_{L^p_{t}L^q_x(\mathbb{T}^{d+1})}\ls_\eps N^{\frac2d+\eps}\|f\|_{L^2(\mathbb{T}^d)},\ \ \alpha\ge2.
\end{equation}
 For $d>4$, this gives the best estimate  up to now, even in the case $\alpha=2$. See e.g. \cite{djlm23} for recent related work on the mixed norm $l^2$ decoupling inequality. It seems reasonable to expect that
\begin{equation}\label{keeltaoconj}	\|e^{it\Delta^{\alpha/2}}f\|_{L^p_{t}L^q_x(\mathbb{T}^{d+1})}\ls_\eps N^{\eps}\|f\|_{L^2(\mathbb{T}^d)},\ \ \alpha\ge2.
\end{equation}
 See \cite[Conjecture 2.6]{MR3374964} for the related Discrete Restriction Conjecture on the flat torus.

 In the following, we use Theorem \ref{fracstr2} to improve Theorem \ref{thm1} on the flat torus.  By Theorem \ref{fracstr2} we can obtain  sharp Strichartz estimates for orthonormal  systems in the subcritical regime. 
\begin{cor}\label{torussub} Let $d\ge 1,\,\alpha \ge2,\ N\ge10$. Suppose $2\le q\le\frac{2(d+2)}{d}$ and $\frac1p=\frac d2(\frac12-\frac1q)$. 
	 Then 
\begin{equation}\label{str-diag}
\Big\|\sum_j  \nu_j |e^{it\Delta^{\alpha/2}} f_j|^2\Big\|_{L^{p/2}_{t}L_x^{q/2}(\mathbb T^{d+1})} \lesssim 
 N^{\sigma} \|\nu\|_{l^\beta} 
 \end{equation}
 holds for  all orthonormal systems $(f_j)_j\subset L^2(\mathbb T^d)$ with $supp\, \hat f_j\subset [-N,N]^d$,  and all sequences $\nu=(\nu_j)_j\in l^\beta$, and all $\sigma \in(0, d]$ and $\beta< \frac{d}{d-\sigma}$.
\end{cor}
The range of $\beta$  is essentially sharp by the necessary condition \eqref{nec}. A remarkable feature  is that  the range $\sigma\in (0,d]$ greatly improves  the one in Theorem \ref{thm2}, and it is essentially optimal by observing the universal bound \eqref{tb}. When $\frac{2(d+2)}d< q<\frac{2(d+1)}{d-1}$, we can also obtain \eqref{str-diag} for all $\sigma\in (2\sigma_1,d]$ and certain $\beta$ by interpolation. Moreover, Nakamura \cite[Theorem 1.4]{MR4068269} proved a similar result  in the case $\alpha=2$ for $p=q=\frac{2(d+2)}d$ and $\sigma\in (0,\frac d{d+2}]$ by directly applying Bourgain-Demeter's $l^2$ decoupling theorem. This result is covered by  Corollary \ref{torussub}.

We may improve the range of $\beta$ in the supercritical regime in Theorem \ref{thm1} on the flat torus, by using the refined estimate at the Keel-Tao endpoint \eqref{keeltao}.

\begin{cor}\label{torussup}
Let $d\ge5$, $\alpha\ge2$ and $N\ge10$. Suppose $\frac{2(d+1)}{d-1}< q\le \frac{2d}{d-2}$ with $\frac1p=\frac d2(\frac12-\frac1q)$. Then \begin{equation}\label{thm6eq}
	\Big\|\sum_j  \nu_j |e^{it\Delta^{\alpha/2}} f_j|^2\Big\|_{L_t^{p/2}L_x^{q/2}(\mathbb T^{d+1})} \lesssim 
	N^{\frac2p}\|\nu\|_{l^\beta}\end{equation}
holds for all orthonormal systems $(f_j)_j\subset L^2(\mathbb T^d)$ with $supp\, \hat f_j\subset [-N,N]^d$, and all sequences $\nu=(\nu_j)_j\in l^\beta$, and all $\beta<\frac{pd(d-3)}{8+pd(d-4)}$.
\end{cor}
Note that $\frac{pd(d-3)}{8+pd(d-4)}>p/2$ whenever $p<\frac{2(d+1)}d$, which is equivalent to $\frac{2(d+1)}{d-1}< q$.  So the range of $\beta$ is larger than the one in the supercritical regime in Theorem \ref{thm1}. Furthermore, if the conjecture \eqref{keeltaoconj} holds, then the conjectural range of $\beta$ in Corollary \ref{torussup} should be $\beta<\frac{d+1}d$, which is exactly the same as the one at the critical point in Theorem \ref{thm1}. See Figure \ref{fig3}.

\subsection{Wave and Klein-Gordon equations}
It is well-known that  the Klein-Gordon propagator $e^{it\sqrt{m^2+\Delta}}$ behaves like the wave propagator $e^{it\sqrt{\Delta}}$ in the high-frequency regime, while it behaves like the Schr\"odinger propagator $e^{it\Delta}$ in the low-frequency regime. For the Strichartz estimates on compact manifolds, the high-frequency regime is more significant. As in the forthcoming Lemma \ref{dislem}, the wave and the Klein-Gordon equations share the same dispersive property on compact manifolds. Indeed, they can be handled in a unified way as pseudodifferential operators, see  Sogge \cite[Chapter 4]{fio}.  So we shall only consider the sharp wave admissible pairs $(p,q)$. By the finite propagation speed property, the Strichartz estimates on compact manifolds are essentially the same as in the Euclidean space.

\begin{figure}[h]
	\centering
	\includegraphics[width=0.7\textwidth]{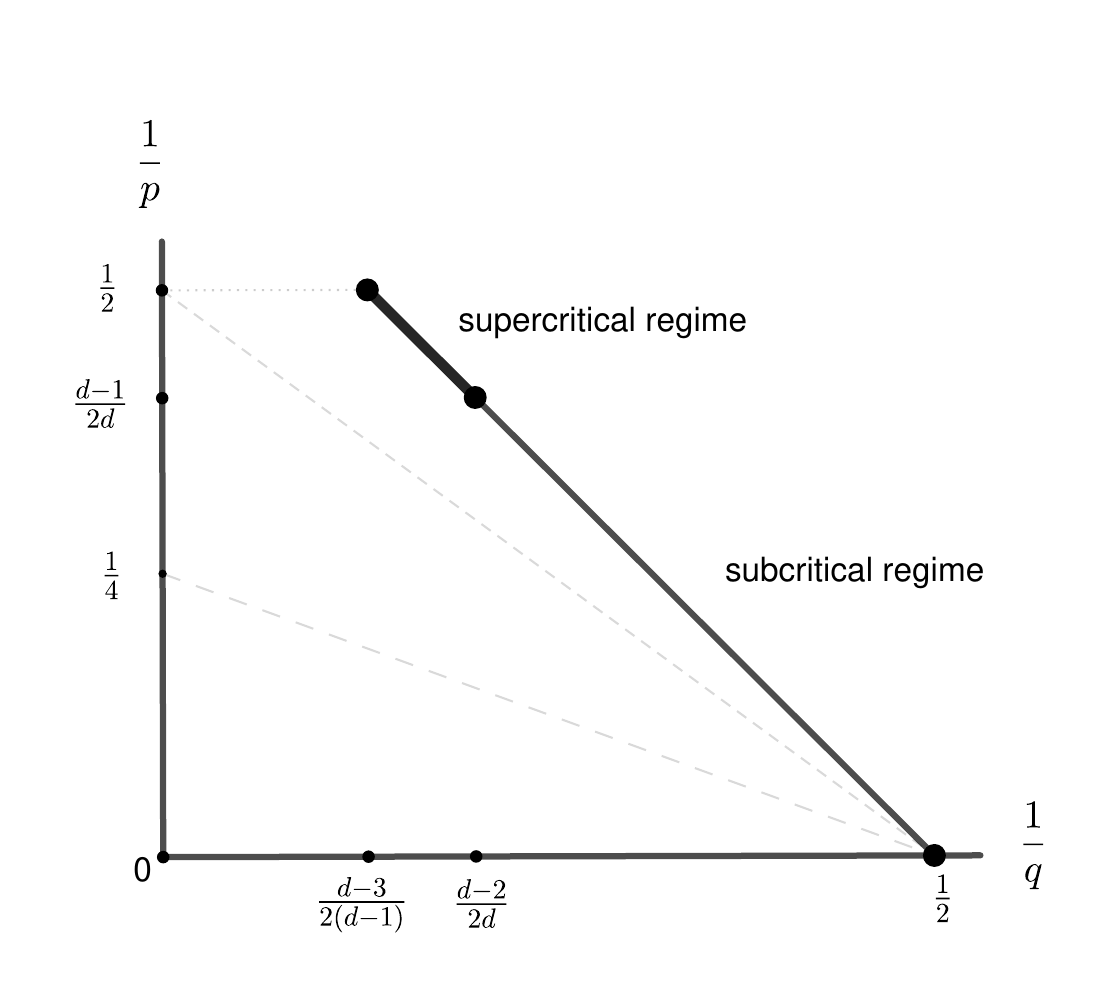}
	\caption{Sharp wave admissible pairs}
	\label{fig2}
\end{figure}

Let $d\ge2$. Suppose $p\ge 2,\ q<\infty$ and $\frac1p=\frac {d-1}2(\frac12-\frac1q)$. We divide these  sharp wave admissible pairs $(p,q)$  in the following four groups.  See Figure \ref{fig2}.

(i) Subcritical regime: $d\ge2$, $2\le q<\frac{2d}{d-2}$

(ii) Critical point: $d\ge3$, $q=\frac{2d}{d-2}$

(iii) Supercritical regime: $d=3$, $\frac{2d}{d-2}< q<\infty$ or $d\ge4$, $\frac{2d}{d-2}< q< \frac{2(d-1)}{d-3}$

(iv) Keel-Tao endpoint: $d\ge4$, $q=\frac{2(d-1)}{d-3}$.

\begin{thm}\label{wkgthm}
	Let $d\ge 2,\ m\ge0,\ N\ge10$. Suppose $p\ge 2,\ q<\infty$ and $\frac1p=\frac {d-1}2(\frac12-\frac1q)$. 
	Let $\sigma_0=\frac2p\frac{d+1}{d-1}$. Then 
	\begin{equation}\label{wkgeq}
		\Big\|\sum_j  \nu_j |e^{it\sqrt{m^2+\Delta}} f_j|^2\Big\|_{L_t^{p/2}L_x^{q/2}(I\times M)} \lesssim N^{\sigma_0}
		\|\nu\|_{l^\beta}\end{equation}
	holds for all orthonormal systems $(f_j)_j$ in $L^2(M)$ with $supp\, \hat f_j\subset \{k:\la_k\le N\}$, and all sequences $\nu=(\nu_j)_j\in l^\beta$, and the following $\beta$ with respect to the pairs $(p,q)$  in the four groups:
	
	{\rm (i)} Subcritical regime: $\beta\le \frac{d-1}{d-1-2/p}=\frac{2q}{q+2}$
	
	{\rm(ii)} Critical point: $\beta<p/2$
	
	{\rm(iii)} Supercritical regime:  $\beta< p/2$

	{\rm(iv)} Keel-Tao endpoint: $\beta=1$.
\end{thm}

The ranges of $\beta$ in (ii)(iii)(iv) are  sharp on any compact manifold by the necessary condition \eqref{nec222}. It is open to show the sharpness for the subcritical regime. Our necessary conditions \eqref{nec} and \eqref{nec222} show the sharpness in the extreme cases $q=2$ and $q=\frac{2d}{d-2}$, so we expect that the intermediate cases are also sharp in some sense. The difficulty to show the sharpness for the subcritical regime also appears in the Euclidean version, see Bez-Lee-Nakamura  \cite[Section 6]{BLN21}.  As before, the optimality of $\beta$ only makes sense when the exponent of $N$  is fixed, so we fix $\sigma_0/2$ to be the Sobolev exponent in the classical (single-function) Strichartz estimates for the wave equation on compact manifolds by Kapitanski \cite{kap89} (see also \cite{cdm23}).

\subsection{Frequency global estimates}By the vector-valued version of the Littlewood-Paley inequality (see \cite[Lemma 1]{sabin16}), one can upgrade the frequency localized estimate  \eqref{eq0} to the frequency global one: for all $s>\sigma/2$,
\begin{equation}\label{eq1}
	\Big\|\sum_j  \nu_j |e^{itP} f_j|^2\Big\|_{L_t^{p/2}L_x^{q/2}(I\times M)} \lesssim
	\|\nu\|_{l^\beta}\end{equation}
holds for all orthonormal systems $(f_j)_j$ in $H^{s}(M)$ and all sequences $\nu=(\nu_j)_j\in l^\beta$. Unlike the single-function case, it seems difficult to take $s=\sigma/2$ in the case of orthonormal data by  Littlewood-Paley theory. Similar difficulties also appear in the Euclidean case. Indeed, Bez-Hong-Lee-Nakamura-Sawano \cite{MR3985036}  observed a crucial fact that on a certain critical line the desired estimates without the frequency localization are not true, see \cite[Prop. 5.2]{MR3985036}. Bez-Lee-Nakamura  \cite{BLN21} achieved frequency global estimates by establishing delicate weighted oscillatory integral estimates and the global dispersive estimates in $\mathbb{R}^d$, see \cite[Prop. 6]{BLN21}. Nevertheless, it seems difficult to proceed in this way on compact manifolds.

\subsection{Organization of the paper} In Section 2, we prove Theorem \ref{thm1} and Theorem \ref{wkgthm} by combining the frequency localized  dispersive estimates for small time intervals with the duality principle due to Frank-Sabin. In Section 3, we prove Strichartz estimates  with variable exponents of $N$ and investigate how  the optimal range of $\beta$ depends on the exponent of $N$. In Section 4, we obtain the improvements on the flat torus by establishing new  decoupling inequalities for certain non-smooth hypersurfaces. In Section 5, we show the optimality of our main theorems. In Section 6, we discuss the applications to the well-posedness of infinite systems of dispersive equations with Hartree-type nonlinearity.

\section{Proof of Strichartz estimates on general manifolds}\label{s-str}
To prove Theorem \ref{thm1} and Theorem \ref{wkgthm}, we shall use the duality principle in Frank-Sabin \cite{MR3730931} that can transfer the orthonormal inequalities to Schatten norm estimates, and the frequency localized dispersive estimates  in Burq-G\'erard-Tzvetkov  \cite[Lemma 2.5]{MR2058384}, Dinh \cite[(3.8)]{dinh17}, Cacciafesta-Danesi-Meng \cite[Prop. 3]{cdm23}, Sogge \cite[Chapter 4]{fio}.

 We recall the definition of the Schatten norm. For $\beta\in [1,\infty)$, $\mathfrak S^\beta=\mathfrak S^\beta(L^2(M))$ denotes the Schatten space based on $L^2(M)$ that is the space of all compact operators $T$ on $L^2(M)$ such that $\text{Tr}|T|^\beta<\infty$ with $|T|=\sqrt{TT^*}$, and its norm is defined by $\|T\|_{\mathfrak S^\beta}=(\text{Tr}|T|^\beta)^{\frac1\beta}$. If $\beta=\infty$, we define $$\|T\|_{\mathfrak S^\infty}=\|T\|_{L^2\to L^2}.$$ Also, $\mathfrak S^2$ is the Hilbert-Schimdt class and the $\mathfrak S^2$ norm is given by  $$\|T\|_{\mathfrak S^2}=\|K\|_{L^2(M\times M)}$$ if $K$ is the integral kernel of $T$. 	See Simon \cite{Simon} for more details on the Schatten classes.

Next, we recall the duality principle in Frank-Sabin \cite[Lemma 3]{MR3730931}.
\begin{lem}[Duality principle]\label{l-dua} Let $p,q,\beta \ge1$. Suppose that $T$ is a bounded operator from $L^2$ to $L_t^pL_x^q$. Then
\[	\|\sum_j\nu_j|Tf_j|^2\|_{L_t^{p/2}L_x^{q/2}}\le C \|\nu\|_{l^\beta}\]
	holds for all orthonormal system $(f_j)_j$ in $L^2$ and all $\nu=(\nu_j)_j$ in $l^\beta$ if and only if
\[	\|WTT^*\overline{W}\|_{\mathfrak S^{\beta'}}\le C \|W\|_{L_t^{2(p/2)'}L_x^{2(q/2)'}}^2\]
	holds for all $W\in L_t^{2(p/2)'}L_x^{2(q/2)'}$.
\end{lem}

\subsection{Proof of Theorem \ref{thm1}}
Let  $\psi\in C_0^\infty(\mathbb{R})$ satisfy 
\[\1_{[-1,1]^d}\le \psi\le \1_{[-2,2]^d}.\]
Let $(e_k)_k$ be an orthonormal eigenbasis in $L^2(M)$ associated with the eigenvalues $(\la_k)_k$ of $\sqrt{\Delta}$. For $f\in L^2(M)$, let
\[\mathcal{E}_N f(t,x)=\sum_{k} \psi(\la_k/N) \hat f(k) e^{ i t\la_k^\alpha}e_k(x).\]
  To prove \eqref{str}, it suffices to show 
\[\Big\|\sum_j  \nu_j |\mathcal E_N f_j|^2\Big\|_{L^{p/2}_tL^{q/2}_x(I\times M)} \lesssim N^{\sigma_0} \|\nu\|_{l^\beta}.\]


At  the endpoint $(p,q)=(\infty,2)$ we have $\beta=1$, and by Minkowski inequality and  Plancherel theorem
\begin{align}\label{end0}
	\|	\sum_{j}\nu_j|\mathcal{E}_N f_j|^2\|_{L_t^{\infty}L_x^{1}}&\le \|\nu\|_{l^1}\sup_j\|\mathcal{E}_N f_j\|_{L_t^\infty L_x^2}^2\ls \|\nu\|_{l^1}.
\end{align}
When $d\ge3$, at the Keel-Tao endpoint $(p,q)=(2,\frac{2d}{d-2})$, we also have $\beta=1$. By the Strichartz estimates \cite[Theorem 1.2]{dinh17} and  Minkowski inequality we have 
\begin{align}
	\|	\sum_{j}\nu_j|\mathcal{E}_N f_j|^2\|_{L_t^{1}L_x^{\frac d{d-2}}}&\le \|\nu\|_{l^1}\sup_j\|\mathcal{E}_N f_j\|_{L_t^2L_x^{\frac{2d}{d-2}}}^2\\\nonumber
	&\ls \|\nu\|_{l^1}\cdot\begin{cases}
		N,\ \ \ \ \ \ \ \ \alpha>1\\
		N^{2-\alpha},\ \ \alpha\in(0,1).
	\end{cases}
\end{align}
 Fix $r\in (d+1,d+2)$. By interpolation, we only need to prove  
\begin{equation}\label{crit}
	\|	\sum_{j}\nu_j|\mathcal{E}_N f_j|^2\|_{L_t^{p/2}L_x^{q/2}}\ls \|\nu\|_{l^\beta}\cdot\begin{cases}
		N^\frac2p,\ \ \ \ \ \ \ \ \alpha>1\\
		N^{\frac{2(2-\alpha)}p},\ \ \alpha\in(0,1).
	\end{cases}
\end{equation}for $\beta'=2(q/2)'=dp/2=r$. The case $d=1,2$ can also be handled similarly. These $(p,q)$ can be close to the critical point $(\frac{2(d+1)}d,\frac{2(d+1)}{d-1})$ as $r\to d+1$.   There is a correspondence between $r$ and the range of $(p,q)$: the subcritical regime ($r>d+1$),  the critical point ($r=d+1$), the supercritical regime ($d<r<d+1$), the Keel-Tao endpoint ($r=d$). 

 Let $\sum_{\ell\ge0}\varphi_\ell(s)=1$ be the Littlewood-Paley decomposition, where $\varphi_\ell\in C_0^\infty$ and for each $\ell>0$, $\varphi_\ell$ is supported in $\{|s|\approx 2^\ell\}$.  Let
\begin{equation}\label{schNL}
	\mathcal{E}_{N,\ell}f(x)=\sum_{k} \psi(\la_k/N)\varphi_\ell(\la_k) e^{i t\la_k^\alpha}\hat f(k)e_k(x).
\end{equation}
When $\alpha>1$, we further split the interval $I$ into $\approx 2^{(\alpha-1)\ell}$ short intervals $\{I_{\ell,n}\}_n$ of length $2^{(1-\alpha)\ell}$. When $\alpha\in (0,1)$, we do not need to split the interval $I$, as it will be clear from the forthcoming Lemma \ref{din}. Then by Minkowski inequality we have
\begin{align}\label{decomp}
\|	\sum_{j}\nu_j|\mathcal{E}_N f_j|^2\|_{L_t^{p/2}L_x^{q/2}}&\le \Big(\sum_{2^\ell\ls N}\|\sum_j\nu_j|\mathcal{E}_{N,\ell}f_j|^2\|_{L^{p/2}_tL_x^{q/2}}^{1/2}\Big)^2\\ \nonumber
&\le\Big(\sum_{2^\ell\ls N}2^{\frac1p(\alpha-1)\ell}\max_n\|\sum_j\nu_j|\mathcal{E}_{N,\ell}f_j|^2\|_{L_t^{p/2}(I_{\ell,n})L_x^{q/2}}^{1/2}\Big)^2.
\end{align}

For any small $\eps>0$, we define an analytic family of  operators $T^\eps_{z,\ell}$ on the strip $$\{z\in\mathbb{C}: -r/2\le \text{Re}z\le 0\}$$ with
the kernels
\[K^\eps_{z,\ell}(t,x,s,y)=\1_{|t|\le 2^{(1-\alpha)\ell}}\1_{|s|\le 2^{(1-\alpha)\ell}}\1_{\eps< |t-s|}(t-s)^{-1-z}\sum_{k} \psi(\la_k/N)^2\varphi_\ell(\la_k)^2 e^{i (t-s)\la_k^\alpha}e_k(x)\overline{e_k(y)}.  \]
By the duality principle in Lemma \ref{l-dua}, we need to  estimate $\|WT^\eps_{-1,\ell}  \overline{W}\|_{\mathfrak S^{r}}$, which follows from the Stein interpolation between the bounds of $\|WT^\eps_{z_1,\ell}  \overline{W}\|_{\mathfrak S^{2}}$ and $\|WT^\eps_{z_2,\ell}  \overline{W}\|_{\mathfrak S^{\infty}}$, where  $z_1=-\frac r2+ib$ and $z_2=ib$ with $b\in\mathbb{R}$.
  
We shall use the frequency localized dispersive estimates  in Burq-G\'erard-Tzvetkov  \cite[Lemma 2.5]{MR2058384} and Dinh \cite[(3.8)]{dinh17}.
\begin{lem}\label{din}
	Let $\alpha\in (0,\infty)\setminus\{1\}$. Let $\varphi\in C_0^\infty(\mathbb{R}\setminus\{0\})$. There exists $t_0>0$ and $C>0$ such that for  any $h\in (0,1]$
\[	\|e^{it\Delta^{\alpha/2}}\varphi(h\sqrt\Delta)f\|_{L^\infty(M)}\le Ch^{-d}(1+|t|h^{-\alpha})^{-d/2}\|f\|_{L^1(M)}\]
	for each $t\in [-t_0h^{\alpha-1},t_0h^{\alpha-1}]$.
\end{lem}

By Lemma \ref{din} we obtain for $|t-s|\ls 2^{(1-\alpha)\ell}$ 
\[	|K^\eps_{z_1,\ell}|\ls |t-s|^{\frac {r-d-2}2}2^{(2-\alpha)d\ell/2}.\]
Then
\begin{align*}
\|WT^\eps_{z_1,\ell}  \overline{W}\|_{\mathfrak S^{2}}&=\|W(t,x)K^\eps_{z_1,\ell}(t,x,s,y)\overline{W}(s,y)\|_{L^2_{t,x,s,y}}\\
&\le C \|W\|_{L_t^{\frac4{r-d}}L_x^2}^22^{(2-\alpha)d\ell/2},
\end{align*} 
where we use the Hardy-Littlewood-Sobolev inequality to estimate the $L_t^{\frac2{r-d}}\to L_t^{(\frac2{r-d})'}$ norm of the convolution operator with the kernel  $|t|^{r-d-2}$. So we require $1<\frac2{r-d}<2$, namely $r\in (d+1,d+2)$. The constant $C$ is independent of $\eps$ and $b$.

Next, by Plancherel theorem we have
\begin{align*}
\|WT^\eps_{z_2,\ell}  \overline{W}\|_{\mathfrak S^{\infty}}&= \|WT^\eps_{z_2,\ell}  \overline{W}\|_{L_{t,x}^2\to L_{t,x}^2}\\
&\le \|W\|_{L^\infty_{t,x}}^2\|T^\eps_{z_2,\ell}\|_{L_{t,x}^2\to L_{t,x}^2}\\
&\le C (1+|b|) \|W\|_{L^\infty_{t,x}}^2
\end{align*} 
where we use the uniform $L^2\rightarrow L^2$ boundedness of the truncated Hilbert transform 
\begin{equation}\label{Hilb}H_b^\eps f(t)=\int_{\eps<|t-s|}\frac{f(s)}{(t-s)^{1-ib}}ds\end{equation}
with
\[	\|H_b^\eps\|_{L^2\to L^2}\le C( 1+|b|)\]
where the constant $C$ is independent of $\eps$ and $b$.
See e.g.  Grafakos \cite[Theorem 5.4.1]{grafbook1}, Vega \cite[p. 204]{vega92}.
Using the Stein interpolation with $\theta=\frac2r$, we get 
\begin{equation}\label{T-1}
	\|WT^\eps_{-1,\ell} \overline{W}\|_{\mathfrak S^{r}}\le C  2^{\frac2p(2-\alpha)\ell} \|W\|_{L_t^{\frac{2r}{r-d}}L_x^{r}(M)}^2 ,
\end{equation}
 since
\[-1=(1-\theta)0+\theta(-\frac r2)\,, \quad \frac1r=\frac{1-\theta}\infty+\frac\theta2\,,\quad \frac1{\frac{2r}{r-d}}=\frac{1-\theta}\infty+\frac{\theta}{\frac4{r-d}}.\]  
The constant $C$ is independent of $\eps$. Let $\eps\to0$ in \eqref{T-1}. Then  by the duality principle in Lemma \ref{l-dua},  we have for $\alpha>1$
\begin{equation}\label{small}
	\|\sum_j\nu_j|\mathcal{E}_{N,\ell}f_j|^2\|_{L_t^{p/2}(I_{\ell,n})L_x^{q/2}}\ls 2^{\frac2p(2-\alpha)\ell}\|\nu\|_{l^\beta},\ \ \forall n.
\end{equation}
Plugging this into \eqref{decomp}, we get \eqref{crit} for $\alpha>1$. The case $\alpha\in (0,1)$ is obtained by directly summing the estimate in \eqref{small} over $\ell$, since we do not need to split the interval $I$.

\subsection{Proof of Theorem \ref{wkgthm}}The proof of Theorem \ref{wkgthm} is similar to that of  Theorem \ref{thm1}. The main difference is that the frequency localized dispersive estimates of the wave and the Klein-Gordon equations in Lemma \ref{dislem} hold for all $|t|\ls1$, which is much better than that of  the fractional Schr\"odinger equations in Lemma \ref{din}. So there is no need to decompose the time interval and the Strichartz estimate has no loss of derivatives  compared to the Euclidean version. Let  $\psi\in C_0^\infty(\mathbb{R})$ satisfy 
\[\1_{[-1,1]^d}\le \psi\le \1_{[-2,2]^d}.\]
Let $(e_k)_k$ be an orthonormal eigenbasis in $L^2(M)$ associated with the eigenvalues $(\la_k)_k$ of $\sqrt{\Delta}$. For $f\in L^2(M)$, let
\[\mathcal{E}_N f(t,x)=\sum_{k} \psi(\la_k/N) \hat f(k) e^{ i t\sqrt{m^2+\la_k^2}}e_k(x).\]
To prove \eqref{wkgeq}, it suffices to show 
\[	\Big\|\sum_j  \nu_j |\mathcal E_N f_j|^2\Big\|_{L^{p/2}_tL^{q/2}_x(I\times M)} \lesssim N^{\frac2p\frac{d+1}{d-1}} \|\nu\|_{l^\beta}.\]


At  the endpoint $(p,q)=(\infty,2)$ we have $\beta=1$, and 
\begin{align*}\label{end1}
	\|	\sum_{j}\nu_j|\mathcal{E}_N f_j|^2\|_{L_t^{\infty}L_x^{1}}&\le \|\nu\|_{l^1}\sup_j\|\mathcal{E}_N f_j\|_{L_x^2}^2\ls \|\nu\|_{l^1}.
\end{align*}
When $d\ge4$, at the Keel-Tao endpoint $(p,q)=(2,\frac{2(d-1)}{d-3})$, we also have $\beta=1$. By the Strichartz estimates \cite[Theorem 1]{cdm23} for single functions, we have
\begin{align*}
	\|	\sum_{j}\nu_j|\mathcal{E}_N f_j|^2\|_{L_t^{1}L_x^{\frac {d-1}{d-3}}}&\le \|\nu\|_{l^1}\sup_j\|\mathcal{E}_N f_j\|_{L_t^2L_x^{\frac{2(d-1)}{d-3}}}^2\\\nonumber
	&\ls N^{\frac{d+1}{d-1}}\|\nu\|_{l^1}.
\end{align*}
Fix $r\in (d,d+1)$. By interpolation, we only need to prove  
\begin{equation}\label{crit2}
	\|	\sum_{j}\nu_j|\mathcal{E}_N f_j|^2\|_{L_t^{p/2}L_x^{q/2}}\ls N^{\frac2p\frac{d+1}{d-1}}\|\nu\|_{l^\beta}
\end{equation}for $\beta'=2(q/2)'=(d-1)p/2=r$. These $(p,q)$ can be close to $(\frac{2d}{d-1},\frac{2d}{d-2})$ as $r\to d$. The case $d=2,3$ can also be handled similarly. There is a correspondence between $r$ and the range of $(p,q)$: the subcritical regime ($r>d$),  the critical point ($r=d$), the supercritical regime ($d-1<r<d$), the Keel-Tao endpoint ($r=d-1$). 

Let $\sum_{\ell\ge0}\varphi_\ell(s)=1$ be the Littlewood-Paley decomposition, where $\varphi_\ell\in C_0^\infty$ and for each $\ell>0$, $\varphi_\ell$ is supported in $\{|s|\approx 2^\ell\}$.  Let
\begin{equation}\label{Enl}
	\mathcal{E}_{N,\ell}f=\sum_{k} \psi(\la_k/N)\varphi_\ell(\la_k)\hat f(k) e^{i t\sqrt{m^2+\la_k^2}}e_k(x).
\end{equation}
For any small $\eps>0$, we define an analytic family of operators $T^\eps_{z,\ell}$ on the strip $$\{z\in\mathbb{C}: -r/2\le \text{Re}z\le 0\}$$ 
with
the kernels
\[K^\eps_{z,\ell}(t,x,s,y)=\1_{\eps< |t-s|}(t-s)^{-1-z}\sum_{k} \psi(\la_k/N)^2\varphi_\ell(\la_k)^2 e^{i (t-s)\sqrt{m^2+\la_k^2}}e_k(x)\overline{e_k(y)}.  \]
By the duality principle in Lemma \ref{l-dua}, we need to estimate $\|WT^\eps_{-1,\ell}  \overline{W}\|_{\mathfrak S^{r}}$ by the Stein interpolation between the bounds of $\|WT^\eps_{z_1,\ell}  \overline{W}\|_{\mathfrak S^{2}}$ and $\|WT^\eps_{z_2,\ell}  \overline{W}\|_{\mathfrak S^{\infty}}$, where  $z_1=-\frac r2+ib$ and $z_2=ib$ with $b\in\mathbb{R}$.

We shall use the frequency localized dispersive estimates for the wave and the Klein-Gordon equations on compact manifolds.  See e.g. Cacciafesta-Danesi-Meng \cite[Prop. 3]{cdm23}, Sogge \cite[Chapter 4]{fio}. 
	\begin{lem}\label{dislem}
		Let $\varphi\in C_0^\infty(\mathbb{R}\setminus\{0\})$. There exists $t_0>0$ and $C>0$ such that for  any $h\in (0,1]$
		\[	\|e^{it\sqrt{m^2+\Delta}}\varphi(h\sqrt\Delta)f\|_{L^\infty(M)}\le Ch^{-d}(1+|t|/h)^{-(d-1)/2}\|f\|_{L^1(M)}\]
		for each $t\in [-t_0,t_0]$.
	\end{lem}
	
	By Lemma \ref{dislem} we obtain for $|t-s|\ls 1$ 
	\[|K^\eps_{z_1,\ell}|\ls |t-s|^{\frac {r-d-1}2}2^{(d+1)\ell/2}.\]
	Then
	\begin{align*}
		\|WT^\eps_{z_1,\ell}  \overline{W}\|_{\mathfrak S^{2}}&=\|W(t,x)K^\eps_{z_1,\ell}(t,x,s,y)\overline{W}(s,y)\|_{L^2_{t,x,s,y}}\\
		&\le C \|W\|_{L_t^{\frac4{r-d+1}}L_x^2}^22^{(d+1)\ell/2},
	\end{align*} 
	where we use the Hardy-Littlewood-Sobolev inequality to estimate the $L_t^{\frac2{r-d+1}}\to L_t^{(\frac2{r-d+1})'}$ norm of the convolution operator with the kernel  $|t|^{r-d-1}$. So we require $1<\frac2{r-d+1}<2$, namely $r\in (d,d+1)$. The constant $C$ is independent of $\eps$ and $b$.

	Next, by Plancherel theorem and the uniform $L^2\rightarrow L^2$ boundedness of the truncated Hilbert transform \eqref{Hilb} we have
	\begin{align*}
		\|WT^\eps_{z_2,\ell}  \overline{W}\|_{\mathfrak S^{\infty}}&= \|WT^\eps_{z_2,\ell}  \overline{W}\|_{L_{t,x}^2\to L_{t,x}^2}\\
		&\le \|W\|_{L^\infty_{t,x}}^2\|T^\eps_{z_2,\ell}\|_{L_{t,x}^2\to L_{t,x}^2}\\
		&\le C(1+|b|) \|W\|_{L^\infty_{t,x}}^2.
	\end{align*} 
	The constant $C$ is independent of $\eps$ and $b$.
	Using the Stein interpolation  we get 
	\begin{equation}\label{T-12}
		\|WT^\eps_{-1,\ell} \overline{W}\|_{\mathfrak S^{r}}\le C  2^{\frac2p\frac{d+1}{d-1}\ell} \|W\|_{L_t^{\frac{2r}{r-d+1}}L_x^{r}(M)}^2.
	\end{equation}
	The constant $C$ is independent of $\eps$. Let $\eps\to 0$ in \eqref{T-12}. Then by the duality principle in Lemma \ref{l-dua}, we have
	\[\|\sum_j\nu_j|\mathcal{E}_{N,\ell}f_j|^2\|_{L_t^{p/2}L_x^{q/2}}\ls 2^{\frac2p\frac{d+1}{d-1}\ell}\|\nu\|_{l^\beta}.\]
	So summing over $2^\ell\ls N$ we get \eqref{crit2}.

\section{Strichartz estimates with variable exponents}
We extend Theorem \ref{thm1} and Theorem \ref{wkgthm} for variable exponents of $N$ and investigate how  the optimal range of $\beta$ depends on the exponent of $N$.
 \subsection{Fractional Schr\"odinger equations} \begin{thm}\label{thm2}
	Let $d\ge 1,\,\alpha\in (0,\infty)\setminus\{1\},\ N\ge10$. Suppose $p\ge 2,\ q<\infty$ and $\frac1p=\frac d2(\frac12-\frac1q)$. Let $\sigma_0$ be defined in \eqref{gam0}.
	Then 
\[	\Big\|\sum_j  \nu_j |e^{it\Delta^{\alpha/2}} f_j|^2\Big\|_{L_t^{p/2}L_x^{q/2}(I\times M)} \lesssim N^{\sigma}
\|\nu\|_{l^\beta}\]
	holds for all orthonormal systems $(f_j)_j$ in $L^2(M)$ with $supp\, \hat f_j\subset \{k:\la_k\le N\}$, and all sequences $\nu=(\nu_j)_j\in l^\beta$, and the following $\sigma$ and $\beta$ with respect to the pairs $(p,q)$  in the four groups:
	
	(i) Subcritical regime: $\sigma\in [\sigma_0,d]$, $\beta\le \frac{d-\sigma_0}{d-\sigma}\beta_*$
	
	(ii) Critical point: $\sigma\in [\sigma_0,d)$, $\beta<\frac{d-\sigma_0}{d-\sigma}\beta_*$, or $\sigma=d$, $\beta\le \infty$
	
	(iii) Supercritical regime:  $\sigma\in [\sigma_0,\sigma_*)$, $\beta< (\frac2p+\sigma_0-\sigma)^{-1}$, or $\sigma\in [\sigma_*,d]$, $\beta\le \frac{d-\sigma_*}{d-\sigma}\beta_*$

	(iv) Keel-Tao endpoint: $\sigma\in [\sigma_0,\sigma_*)$, $\beta\le (\frac2p+\sigma_0-\sigma)^{-1}$, or $\sigma\in [\sigma_*,d]$, $\beta\le \frac{d-\sigma_*}{d-\sigma}\beta_*$. 
	
	Here \begin{align}
		\beta_*=\frac{d}{d-2/p}=\frac{2q}{q+2},\ \ 	\sigma_*=\sigma_0+\frac2p-\frac1{\beta_*}.
	\end{align}
\end{thm}

When $\alpha>1$,  the ranges of $\beta$ in (i)(ii) are essentially sharp  by the necessary conditions  \eqref{nec}, and the ranges of $\beta $ for $\sigma\in [\sigma_0,\sigma_*]$ in (iii)(iv)  are also essentially sharp by the necessary condition  \eqref{nec22}. It is open to show the ranges of $\beta$ for $\sigma\in [\sigma_*,d]$ in (iii)(iv) are sharp. Since they are sharp in the extreme cases $\sigma=\sigma_*$ and $\sigma=d$ by the necessary condition \eqref{nec22}, we expect that they are also sharp for the intermediate values.

Theorem \ref{thm2} can be deduced from interpolation between Theorem \ref{thm1}, the universal bound \eqref{tb}, and the ``kink point'' estimate 
\begin{equation}\label{kink}
	\|	\sum_{j}\nu_j|e^{it\Delta^{\alpha/2}} f_j|^2\|_{L_t^{p/2}L_x^{q/2}}
	\ls\|\nu\|_{l^{\beta_*}}\cdot\begin{cases}
		N^{\sigma_*},\ \ \ \ \ \ \ \ \ \ q>\frac{2(d+1)}{d-1}\\
		N^{\sigma_*}(\log N),\ q=\frac{2(d+1)}{d-1}.
	\end{cases}
\end{equation}
By the necessary condition \eqref{nec22}, the estimate \eqref{kink} is sharp when $\alpha>1$, in the sense that $\beta_*$ cannot be replaced by any $\beta>\beta_*$. Moreover, Nakamura \cite[Theorem 5.1]{MR4068269} obtained similar estimates for the Schr\"odinger propagator $e^{it\Delta}$ on the flat torus with an $\eps$-loss.

In Theorem \ref{thm2}, the subcritical regime directly follows from interpolation between Theorem \ref{thm1} and the universal bound \eqref{tb}. Nevertheless, this  interpolation argument is not enough to give a sharp range of $\beta$ in the supercritical regime.  Unlike the proof of Theorem \ref{thm1}, we shall further  decompose the kernel of $\mathcal{E}_{N,\ell}\mathcal{E}_{N,\ell}^*$ to establish the kink point estimate \eqref{kink}, where $\mathcal{E}_{N,\ell}$ is given by \eqref{schNL}.

\subsection{Proof of Theorem \ref{thm2}}Let $\sum_{m\in \mathbb{Z}}\tilde \varphi_m(s)=1$ be the Littlewood-Paley decomposition, where each $\tilde \varphi_m\in C_0^\infty$ is supported in $\{|s|\approx 2^m\}$. We define the operator $T_{\ell,m}$ with the kernel 
\[K_{\ell,m}(t,x,s,y)=\1_{|t|\le 2^{(1-\alpha)\ell}}\1_{|s|\le 2^{(1-\alpha)\ell}}\tilde \varphi_m(|t-s|)\sum_{k} \psi(\la_k/N)^2\varphi_\ell(\la_k)^2 e^{i (t-s)\la_k^\alpha}e_k(x)\overline{e_k(y)}.  \]
This comes from the dyadic decomposition of the operator $\mathcal{E}_{N,\ell}\mathcal{E}_{N,\ell}^*$.
Fix $r\in [d,d+1]$. Let $\beta'=2(q/2)'=dp/2=r$.   As before, we shall estimate $\|WT_{\ell,m}  \overline{W}\|_{\mathfrak S^{r}}$ by interpolation between the $\mathfrak S^2$ norm and the $\mathfrak S^\infty$ norm. 

We first estimate the $\mathfrak S^2$ norm. By Lemma \ref{din} we obtain for $ 2^m\ls 2^{(1-\alpha)\ell}$ 
\[	|K_{\ell,m}|\ls \min\{2^{-dm/2}2^{(2-\alpha)d\ell/2},\ 2^{d\ell}\} .\]
Then  we have
\begin{align}\label{s2norm}
	\|WT_{\ell,m}  \overline{W}\|_{\mathfrak S^{2}}&=\|W(t,x)K_{\ell,m}(t,x,s,y)\overline{W}(s,y)\|_{L^2_{t,x,s,y}}\\\nonumber
&\ls \|W\|_{L_t^{\frac4{r-d}}L_x^2}^2 2^{m/2}2^{(\alpha-1)\ell(r-d-1)/2}\cdot \min\{2^{-dm/2}2^{(2-\alpha)d\ell/2},\ 2^{d\ell}\}.
\end{align} 
Here we require $\frac{2}{r-d}\ge2$, namely $d\le r\le d+1$, since we use the inequality 
\begin{equation}\label{conv}
	\Big|\int_I\int_I f(t)h(t-s)g(s)dsdt\Big|\le |I|^{1-\frac2p}\|h\|_{L^1}\|f\|_{L^p}\|g\|_{L^p},\ \ \forall p\ge2.
\end{equation}It simply follows from H\"older inequality.

Next, we estimate the $\mathfrak S^\infty$ norm.   Then  by Plancherel theorem we have
\begin{align*}
	\|WT_{\ell,m} \overline{W}\|_{\mathfrak S^{\infty}}&= \|WT_{\ell,m}  \overline{W}\|_{L_{t,x}^2\to L_{t,x}^2}\\
	&\le \|W\|_{L^\infty_{t,x}}^2\|T_{\ell,m}\|_{L_{t,x}^2\to L_{t,x}^2}\\
	&\ls 2^m\|W\|_{L_{t,x}^\infty}^2.
\end{align*} 
Then  interpolation gives
\begin{align*}
	\|WT_{\ell,m} \overline{W}\|_{\mathfrak S^{r}}\ls \|W\|_{L_t^{\frac{2r}{r-d}}L_x^r}^22^{\frac\ell r(\alpha-1)(r-d-1)}\min\{2^{m(1-\frac{d+1}r)}2^{\frac{d\ell}r(2-\alpha)},2^{m(1-\frac1r)}2^{2d\ell/r}\}.
\end{align*} 
 When $\alpha>1$, summing over $2^m\ls 2^{(1-\alpha)\ell}$ and by the duality principle, we have for each small interval $I_{\ell,n}$
\begin{equation}\label{piece}
	\|\sum_j\nu_j|\mathcal{E}_{N,\ell}f_j|^2\|_{L_t^{p/2}(I_{\ell,n})L_x^{q/2}}\ls\|\nu\|_{l^\beta}\cdot  \begin{cases}2^{\frac \ell r(3d-\alpha d+1-r)},\ \ \ \ \ \ \ r\in [d,d+1)\\
		(1+\ell) 2^{\frac \ell r(3d-\alpha d+1-r)},\ \ r=d+1.
		\end{cases} 
\end{equation}
Plugging this into \eqref{decomp}, we get for $\alpha>1$
\[(\sum_{2^\ell\ls N}2^{\frac{(\alpha-1)d\ell}{2r}}2^{\frac{\ell}{2r}(3d-\alpha d+1-r)})^2\ls N^{\frac{2d+1}r-1}=N^{\frac2p+\frac{2(d+1)}{pd}-1}.\]
When $\alpha\in (0,1)$, we have no need to split the interval  and we directly sum the estimates \eqref{piece} over $\ell$ to obtain
\[\sum_{2^\ell\ls N}2^{\frac{\ell}{r}(3d-\alpha d+1-r)}\ls N^{\frac{(3-\alpha)d+1}r-1}=N^{\frac{2(2-\alpha)}p+\frac{2(d+1)}{pd}-1}.\]
Thus, if we define 
\begin{align*}
\beta_*=\frac{d}{d-2/p},\ \ 	\sigma_*=\sigma_0+\frac{2(d+1)}{pd}-1,
\end{align*}with $\sigma_0$ given in \eqref{gam0}, then we have
\begin{equation}\label{decomp2}
	\|	\sum_{j}\nu_j|e^{it\Delta^{\alpha/2}} f_j|^2\|_{L_t^{p/2}L_x^{q/2}}
	\ls\|\nu\|_{l^{\beta_*}}\cdot\begin{cases}
		N^{\sigma_*},\ \ \ \ \ \ \ \ \ \ r\in [d,d+1)\\
		N^{\sigma_*}(\log N),\ r=d+1.
	\end{cases}
\end{equation}
Note that $\sigma_*=\sigma_0$ is equivalent to $p=\frac{2(d+1)}d$. By interpolation between \eqref{decomp2}, \eqref{tb} and Theorem \ref{thm1}, we obtain Theorem \ref{thm2}.

 In particular, when $d\ge1$, at the critical point $(p,q)=(\frac{2(d+1)}d,\frac{2(d+1)}{d-1})$, we have $r=d+1$, $\beta=\frac{d+1}d$ and
\begin{equation}\label{crib}
		\|	\sum_{j}\nu_j|e^{it\Delta^{\alpha/2}} f_j|^2\|_{L_t^{p/2}L_x^{q/2}}
	\ls \|\nu\|_{l^\beta}\cdot \begin{cases}
		N^{\frac2p}(\log N),\ \ \ \ \ \ \ \alpha>1\\
		N^{\frac{2(2-\alpha)}p}(\log N),\ \ \alpha\in (0,1).
	\end{cases}
\end{equation}
When $d\ge2$, at the Keel-Tao endpoint $(p,q)=(2,\frac{2d}{d-2})$, we have $r=d$, $\beta=\frac{d}{d-1}$ and 
\begin{equation}\label{ktb}
		\|	\sum_{j}\nu_j|e^{it\Delta^{\alpha/2}} f_j|^2\|_{L_t^{p/2}L_x^{q/2}}
	\ls \|\nu\|_{l^\beta}\cdot \begin{cases}
		N^{\frac2p+\frac1d},\ \ \ \ \ \ \ \ \alpha>1\\
		N^{\frac{2(2-\alpha)}p+\frac1d},\ \ \alpha\in (0,1).
	\end{cases}
\end{equation}
Note that when $\alpha>1$ this estimate is sharp by the necessary condition \eqref{nec22}.

$ $

\noindent \textbf{Remark 1.}  To compare this method with the proof of Theorem \ref{thm1}, we remark that one can slightly modify this method to handle the subcritical regime ($r>d+1$) up to a log loss:
\begin{equation}\label{subloss}
	\|	\sum_{j}\nu_j|e^{it\Delta^{\alpha/2}} f_j|^2\|_{L_t^{p/2}L_x^{q/2}}
	\ls \|\nu\|_{l^\beta}\cdot \begin{cases}
		N^{\frac2p}(\log N),\ \ \ \ \ \ \ \ \alpha>1\\
		N^{\frac{2(2-\alpha)}p}(\log N),\ \ \ \alpha\in(0,1).
	\end{cases}
\end{equation}
Indeed, for $r\in (d+1,d+2)$ by Young's inequality, we can replace \eqref{s2norm} by 
\begin{align*}
	\|WT_{\ell,m}  \overline{W}\|_{\mathfrak S^{2}}&=\|W(t,x)K_{\ell,m}(t,x,s,y)\overline{W}(s,y)\|_{L^2_{t,x,s,y}}\\\nonumber
	&\ls \|W\|_{L_t^{\frac4{r-d}}L_x^2}^2 2^{\frac{2+d-r}2m}\cdot \min\{2^{-dm/2}2^{(2-\alpha)d\ell/2},\ 2^{d\ell}\}.
\end{align*} 
Then repeating the interpolation argument above we can obtain \eqref{subloss}. So we expect that this method  gives essentially  sharp bounds for all admissible $(p,q)$.

$ $

\noindent \textbf{Remark 2.} It is natural to ask whether one can remove log factor in \eqref{crib} at the critical point $(p,q)=(\frac{2(d+1)}d,\frac{2(d+1)}{d-1})$. Recall that  it is not removable in the Euclidean version, see \cite{MR3254332, MR3985036, BLN21}. However, surprisingly it can be removed for the Schr\"odinger propagator $e^{it\Delta}$ on the one dimensional flat torus by following the spirit of the Hardy-Littlewood circle method. See  Nakamura \cite[Theorem 1.6]{MR4068269}.
\subsection{Wave and Klein-Gordon equations}
We may also extend Theorem \ref{wkgthm} for variable exponents of $N$.
	 \begin{thm}\label{wkgthm2}
		Let $d\ge 2,\ N\ge10$. Let  $m\ge0$.   Suppose $p\ge 2,\ q<\infty$ and $\frac1p=\frac {d-1}2(\frac12-\frac1q)$. Let $\sigma_0=\frac2p\frac{d+1}{d-1}$.
		Then 
		\[	\Big\|\sum_j  \nu_j |e^{it\sqrt{m^2+\Delta}} f_j|^2\Big\|_{L_t^{p/2}L_x^{q/2}(I\times M)} \lesssim N^{\sigma}
		\|\nu\|_{l^\beta}\]
		holds for all orthonormal systems $(f_j)_j$ in $L^2(M)$ with $supp\, \hat f_j\subset \{k:\la_k\le N\}$, and all sequences $\nu=(\nu_j)_j\in l^\beta$, and the following $\sigma$ and $\beta$ with respect to the pairs $(p,q)$  in the four groups:
		
		(i) Subcritical regime: $\sigma\in [\sigma_0,d]$, $\beta\le \frac{d-\sigma_0}{d-\sigma}\beta_*$
		
		(ii) Critical point: $\sigma\in [\sigma_0,d)$, $\beta<\frac{d-\sigma_0}{d-\sigma}\beta_*$, or $\sigma=d$, $\beta\le \infty$
		
		(iii) Supercritical regime:  $\sigma\in [\sigma_0,\sigma_*)$, $\beta< (\frac2p+\sigma_0-\sigma)^{-1}$, or $\sigma\in [\sigma_*,d]$, $\beta\le \frac{d-\sigma_*}{d-\sigma}\beta_*$

		(iv) Keel-Tao endpoint: $\sigma\in [\sigma_0,\sigma_*)$, $\beta\le (\frac2p+\sigma_0-\sigma)^{-1}$, or $\sigma\in [\sigma_*,d]$, $\beta\le \frac{d-\sigma_*}{d-\sigma}\beta_*$. 
		Here \begin{align}
			\beta_*=\frac{d-1}{d-1-2/p}=\frac{2q}{q+2},\ \ 	\sigma_*=\sigma_0+\frac2p-\frac1{\beta_*}.
		\end{align}
	\end{thm}
	The ranges of $\beta$ in (ii)(iii)(iv) are essentially sharp on the sphere for $\sigma\in [\sigma_0,\sigma_*]$ by the necessary conditions \eqref{nec222}. It is open to show other ranges of $\beta$ are sharp. By the necessary conditions \eqref{nec} and \eqref{nec222}, the range of $\beta$ in (i) is sharp in the extreme cases $q=2$ and $q=\frac{2d}{d-2}$. By the necessary condition \eqref{nec222}, the ranges of $\beta$ for $\sigma\in [\sigma_*,d]$ in (iii)(iv) are sharp in the extreme cases $\sigma=\sigma_*$ and $\sigma=d$. So we expect that they are also sharp for the intermediate values.

	Theorem \ref{wkgthm2} can be deduced from interpolation between Theorem \ref{wkgthm}, the universal bound \eqref{tb}, and the ``kink point'' estimate 
	\begin{equation}\label{kink2}
		\|	\sum_{j}\nu_j|e^{it\sqrt{m^2+\Delta}} f_j|^2\|_{L_t^{p/2}L_x^{q/2}}
		\ls\|\nu\|_{l^{\beta_*}}\cdot\begin{cases}
			N^{\sigma_*},\ \ \ \ \ \ \ \ \ \ q>\frac{2d}{d-2}\\
			N^{\sigma_*}(\log N),\ q=\frac{2d}{d-2}.
		\end{cases}
	\end{equation}
	By the necessary condition \eqref{nec222}, the estimate \eqref{kink2} is sharp in the sense that $\beta_*$ cannot be replaced by any $\beta>\beta_*$.

	The subcritical regime in Theorem \ref{wkgthm2}  directly follows from interpolation between Theorem \ref{wkgthm} and the universal bound \eqref{tb}. To handle the supercritical regime, we shall modify the proof of Theorem \ref{thm2} to establish the kink point estimate \eqref{kink2}. The argument is similar to the proof of \eqref{kink}, which dyadically decomposes the kernel of $\mathcal{E}_{N,\ell}\mathcal{E}_{N,\ell}^*$ with respect to $|t-s|$.
\subsection{Proof of Theorem \ref{wkgthm2}}	
	Let $\sum_{n\in \mathbb{Z}}\tilde \varphi_n(s)=1$ be the Littlewood-Paley decomposition, where each $\tilde \varphi_n\in C_0^\infty$ is supported in $\{|s|\approx 2^n\}$.  We define the operator $T_{\ell,n}$ with the kernel 
	\[K_{\ell,n}(t,x,s,y)=\tilde \varphi_n(|t-s|)\sum_{k} \psi(\la_k/N)^2\varphi_\ell(\la_k)^2 e^{i (t-s)\sqrt{m^2+\la_k^2}}e_k(x)\overline{e_k(y)}.  \]
	This comes from the dyadic decomposition of the operator $\mathcal{E}_{N,\ell}\mathcal{E}_{N,\ell}^*$ with $\mathcal{E}_{N,\ell}$ given by \eqref{Enl}.
	Fix $r\in [d-1,d]$. Let $\beta'=2(q/2)'=(d-1)p/2=r$.   As before, we shall estimate $\|WT_{\ell,n}  \overline{W}\|_{\mathfrak S^{r}}$ by interpolation between the $\mathfrak S^2$ norm and the $\mathfrak S^\infty$ norm. 
	
	We first estimate the $\mathfrak S^2$ norm. By Lemma \ref{dislem} we obtain for $ 2^n\ls 1$ 
	\[	|K_{\ell,n}|\ls \min\{2^{-(d-1)n/2}2^{(d+1)\ell/2},\ 2^{d\ell}\} .\]
	Then  we have
	\begin{align*}
		\|WT_{\ell,n}  \overline{W}\|_{\mathfrak S^{2}}&=\|W(t,x)K_{\ell,n}(t,x,s,y)\overline{W}(s,y)\|_{L^2_{t,x,s,y}}\\\nonumber
		&\ls \|W\|_{L_t^{\frac4{r-d+1}}L_x^2}^2 2^{n/2}\cdot \min\{2^{-(d-1)n/2}2^{(d+1)\ell/2},\ 2^{d\ell}\}.
	\end{align*} 
	Here we require $\frac{2}{r-d+1}\ge2$, namely $d-1\le r\le d$, by using the inequality \eqref{conv}.
	
	Next, we estimate the $\mathfrak S^\infty$ norm.   Then  by Plancherel theorem we have
	\begin{align*}
		\|WT_{\ell,n} \overline{W}\|_{\mathfrak S^{\infty}}&= \|WT_{\ell,n}  \overline{W}\|_{L_{t,x}^2\to L_{t,x}^2}\\
		&\le \|W\|_{L^\infty_{t,x}}^2\|T_{\ell,n}\|_{L_{t,x}^2\to L_{t,x}^2}\\
		&\ls 2^n\|W\|_{L_{t,x}^\infty}^2.
	\end{align*} 
	Then  interpolation gives
	\begin{align*}
		\|WT_{\ell,n} \overline{W}\|_{\mathfrak S^{r}}\ls \|W\|_{L_t^{\frac{2r}{r-d+1}}L_x^r}^22^{n(1-\frac1r)}\min\{2^{-(d-1)n/r}2^{(d+1)\ell/r},\ 2^{2d\ell/r}\}.
	\end{align*} 
	Summing over $2^n\ls 1$ and by the duality principle in Lemma \ref{l-dua}, we have
	\[\|\sum_j\nu_j|\mathcal{E}_{N,\ell}f_j|^2\|_{L_t^{p/2}L_x^{q/2}}\ls\|\nu\|_{l^\beta}\cdot  \begin{cases}2^{(\frac{2d+1}r-1)\ell},\ \ \ \ \ \ \ \ \ \ r\in [d-1,d)\\
		(1+\ell) 2^{(\frac{2d+1}r-1)\ell},\ \ r=d.
	\end{cases} \]
	Recall that $r=\beta'=(d-1)p/2$. Thus, if we define 
	\begin{align*}
		\beta_*=\frac{d-1}{d-1-2/p},\ \ 	\sigma_*=\frac{2(d+1)}{p(d-1)}+\frac{2d}{p(d-1)}-1,
	\end{align*} then summing over $\ell$ we have
	\[	\|	\sum_{j}\nu_j|\mathcal{E}_{N,\ell} f_j|^2\|_{L_t^{p/2}L_x^{q/2}}
	\ls\|\nu\|_{l^{\beta_*}}\cdot\begin{cases}
		N^{\sigma_*},\ \ \ \ \ \ \ \ \ \ r\in [d-1,d)\\
		N^{\sigma_*}(\log N),\ r=d.
	\end{cases}\]
	This proves \eqref{kink2}. Note that $\sigma_*=\frac2p\frac{d+1}{d-1}$ is equivalent to $p=\frac{2d}{d-1}$. By interpolation between \eqref{kink2}, \eqref{tb} and Theorem \ref{wkgthm}, we obtain Theorem \ref{wkgthm2}.
	
	In particular, when $d\ge2$, at the critical point $(p,q)=(\frac{2d}{d-1},\frac{2d}{d-2})$, we have $r=d$, $\beta=\frac{d}{d-1}$ and
	\[\|	\sum_{j}\nu_j|e^{it\sqrt{m^2+\Delta}} f_j|^2\|_{L_t^{p/2}L_x^{q/2}}
	\ls N^{\frac2p\frac{d+1}{d-1}}(\log N)\|\nu\|_{l^\beta}.\]
	When $d\ge3$, at the Keel-Tao endpoint $(p,q)=(2,\frac{2(d-1)}{d-3})$, we have $r=d-1$, $\beta=\frac{d-1}{d-2}$ and 
	\[\|	\sum_{j}\nu_j|e^{it\sqrt{m^2+\Delta}} f_j|^2\|_{L_t^{p/2}L_x^{q/2}}
	\ls N^{\frac2p\frac{d+1}{d-1}+\frac1{d-1}}\|\nu\|_{l^\beta}.\]
	Note that this estimate is sharp by the necessary condition \eqref{nec222}.

\section{Decoupling inequalities and improvements on the flat torus}
Let $d\ge1$. For $c_1,...,c_{d+1}>0$, let $R=[-c_1,c_1]\times...\times[-c_{d+1},c_{d+1}]$ be a rectangular box. We shall use two weight functions associated with $R$
\[\omega_R(x)=(1+\sum_{j=1}^{d+1}|x_j|/c_j)^{-10d}\]
\[\tilde\omega_R(x)=(1+\sum_{j=1}^{d+1}|x_j|/c_j)^{-8d}.\]
Let $\alpha>1$ and $S=\{(y,|y|^\alpha)\in \mathbb{R}^{d+1}:y\in [-1,1]^d\}$. We define the extension operator
\[E_Q g(x)=\int_Q g(y)e(x_1y_1+...+x_dy_d+x_{d+1}|y|^\alpha )dy\]
where $Q$ is some subset in $\mathbb{R}^{d}$ and $e(z)=e^{2\pi i z}$. Let $\text{Part}_{\delta^{1/2}}([-1,1]^d)$ denote a partition of $[-1,1]^d$ into cubes of side length $\delta^{1/2}$. We first prove the following decoupling inequality for the hypersurface $S$.

\begin{thm}\label{dec0}
	Let $\alpha>1$ and $2\le p\le \frac{2(d+2)}{d}$. Then we have for all $\eps>0$,
	\begin{equation}\label{dec}
		\|E_{[-1,1]^d} g\|_{L^p(B_R)}\ls_\eps \delta^{-\eps}\Big(\sum_{\Delta\in \text{Part}_{\delta^{1/2}}([-1,1]^d)}	\|E_{\Delta} g\|_{L^p(\omega_{B_R})}^2\Big)^\frac12
	\end{equation}
	where $B_R$ is a ball of radius $R\ge \delta^{-{\max\{1,\alpha/2\}}}$.
\end{thm}
 The proof extends the strategy of the one dimensional case in \cite{MR4078083}, which used a piece of parabola to locally approximate the curve  and then apply Bourgain-Demeter's decoupling theorem. For recent works on decoupling inequalities for smooth hypersurfaces with vanishing Gaussian curvature, see e.g. Demeter \cite[Section 12.6]{MR3971577}, Yang \cite{Yang21}, Li-Yang  \cite{MR4458402,li2024decouplingsmoothsurfacesmathbbr3} and Guth-Maldague-Oh \cite{guth2024l2decouplingtheoremsurfaces}.  Now we use this inequality to obtain the Strichartz estimates, and postpone its proof to the end of this section. 
 
 Theorem \ref{dec0} immediately implies a discrete restriction estimate.

\begin{cor}\label{corsum}
	Let $\alpha>1$ and $2\le p\le \frac{2(d+2)}{d}$. Let $\Lambda$ be a $N^{-1}$-separated set in $[-1,1]^d$. Then we have for all $\eps>0$,
	\begin{equation}\label{dec2}
		\Big(\frac1{|B_R|}\int_{B_R}\Big|\sum_{\xi\in \Lambda}a_\xi e(x_1\xi_1+...+x_d\xi_d+x_{d+1}|\xi|^\alpha)\Big|^pdx\Big)^\frac1p\ls_\eps N^{\eps}\Big(\sum_{\xi\in \Lambda}|a_\xi|^2\Big)^\frac12
	\end{equation}
	where $B_R$ is a ball of radius $R\ge N^{\max\{2,\alpha\}}$.
\end{cor}
This discrete restriction estimate together with the trivial endpoint $p=\infty$ estimate directly implies the Strichartz estimate on the flat torus in Theorem \ref{fracstr}. 
\begin{thm}\label{fracstr}
	Let $d\ge 1,\,\alpha >1,\ N\ge10$.  Let $f\in L^2(\mathbb{T}^d)$ with $supp \hat f\subset [-N,N]^d$. Then  we have for all $\eps>0$ and $2\le q\le \infty$
\[\|e^{it\Delta^{\alpha/2}}f\|_{L^q_{t,x}(\mathbb{T}^{d+1})}\ls_\eps N^{\sigma_1} \|f\|_{L^2(\mathbb{T}^d)}\]
	where
\begin{equation}\label{sigma11}
	\sigma_1(\alpha)=\begin{cases}
		\max\{(2-\alpha)\frac d2(\frac12-\frac1q),\frac d2-\frac{d+\alpha}q\}+\eps,\ \ 1<\alpha<2\\
		\max\{0,\frac d2-\frac{d+2}q\}+\eps,\ \ \quad\quad\quad\quad\quad\quad\quad\quad\alpha\ge2.
	\end{cases}
\end{equation}
\end{thm}
To our best knowledge, the case $\alpha=2$ is due to Bourgain-Demeter \cite[Theorem 2.4]{MR3374964}. After completing this paper, we learned that results essentially equivalent to our  Theorem \ref{fracstr} were recently obtained  by  Schippa \cite[Prop. 1.1]{sch19}, and the special case $d=1$ and $p=q=4$ was  obtained earlier by Demibras-Erdo\u gan-Tzirakis \cite[Theorem 4]{det16}).  The proof of Theorem \ref{fracstr} is standard. For completeness, we include a short proof here.

$ $

\noindent \textbf{Proof of Theorem \ref{fracstr}.} Let $f(x)=\sum_{k\in\mathbb{Z}^d:|k|\le N}a_ke(k\cdot x)$. Let $\Lambda=\frac1N\mathbb{Z}^d$ and $\xi=k/N$. Let  $I$ be a  bounded interval of length $|I|>0$, By scaling and periodicity in the space variable, we can use Corollary \ref{corsum} with $p=\frac{2(d+2)}d$ to get

\begin{align*}
	\int_{I\times \mathbb{T}^d}|e^{it\Delta^{\alpha/2}}f|^p&\ls
	 \frac1{(N|I|)^{d}N^\alpha}\int_{|t|\ls N^\alpha |I|,\ |y|\ls N|I|}\Big|\sum_{\xi\in\Lambda}a_{N\xi} e(y\cdot \xi+t|\xi|^\alpha)\Big|^pdtdy\\
	 &\ls  \frac1{(N|I|)^{d}N^\alpha}\frac1{N^{d(\alpha-1)}}\int_{|t|\ls N^\alpha |I|,\ |y|\ls N^\alpha|I|}\Big|\sum_{\xi\in\Lambda}a_{N\xi} e(y\cdot \xi+t|\xi|^\alpha)\Big|^pdtdy\\
	 &\ls \frac1{(N|I|)^{d}N^\alpha}\frac1{N^{d(\alpha-1)}}\cdot (N^\alpha |I|)^{d+1}N^\eps\|f\|_{L^2(\mathbb{T}^d)}^p\\
	 &\approx |I|N^\eps\|f\|_{L^2(\mathbb{T}^d)}^p.
\end{align*}
Here we require $N^\alpha |I|\gs N^{\max\{2,\alpha\}}$ in the third  inequality. So we fix $|I|=N^{2-\alpha}$ when $\alpha<2$ and $|I|=1$ when $\alpha\ge2$. Then Theorem \ref{fracstr} follows from interpolation between $p=\frac{2(d+2)}d$ and the  endpoints $p=2,\infty$.

As a consequence, we obtain  estimates for the sharp Schr\"odinger admissible pairs.
\begin{thm}\label{fracstr2}
	Let $d\ge 1,\,\alpha \in(0,\infty)\setminus\{1\},\ N\ge10$.  Let $f\in L^2(\mathbb{T}^d)$ with $supp \hat f\subset [-N,N]^d$. Suppose $p\ge 2,\ q<\infty$ and $\frac1p=\frac d2(\frac12-\frac1q)$.   Then  for all $\eps>0$, 
	\begin{equation}\label{dec4}
		\|e^{it\Delta^{\alpha/2}}f\|_{L^p_{t}L^q_x(\mathbb{T}^{d+1})}\ls_\eps N^{\sigma_2}\|f\|_{L^2(\mathbb{T}^d)}
	\end{equation}
	where 
	\begin{equation}\label{sigma2}
		\sigma_2(\alpha)=\begin{cases}\sigma_0(\alpha)/2,\ \ \ \ \ \ \ \quad\quad\quad 0<\alpha<1\\
			\min\{\sigma_0(\alpha)/2,\ \sigma_1(\alpha)\},\ \ \ \ \ \alpha>1,
		\end{cases}
	\end{equation}
	with $\sigma_0$ and $\sigma_1$ given by \eqref{gam0} and \eqref{sigma11}.
\end{thm}

$ $

\noindent\textbf{Proof of Theorem \ref{fracstr2}.} The exponent $\sigma_0(\alpha)/2$ follows from Dinh \cite[Theorem 1.2]{det16}. We just need to handle $\alpha>1$. We use H\"older inequality to control $L^p_tL^q_x$ norm by $L_t^qL_x^q$ norm in Theorem \ref{fracstr} if $p\le q$. The remaining case $p>q$ follows from interpolation between $p=q=\frac{2(d+2)}d$ and the trivial endpoint $(p,q)=(\infty,2)$. These give the exponent $\sigma_1(\alpha)$.

$ $

\noindent\textbf{Examples, sharpness and improvements.} Consider the function $f(x)=\sum_{k\in\mathbb{Z}^d:|k|\le N}e^{2\pi ik\cdot x}$. We have $\|f\|_{L^2(\mathbb{T}^d)}\approx N^{d/2}$ and
\[\|e^{it\Delta^{\alpha/2}}f\|_{L^p_{t}L^q_x(\mathbb{T}^{d+1})}\ge 	\|e^{it\Delta^{\alpha/2}}f\|_{L^p_{t}(\{|t|<N^{-\alpha}\})L^q_x(\{|x|<N^{-1}\})}\approx N^{d-\frac{d}q-\frac{\alpha}p}.\]
This implies a lower bound \begin{equation}\label{lowb}
	\|e^{it\Delta^{\alpha/2}}f\|_{L^p_{t}L^q_x(\mathbb{T}^{d+1})}\gs N^{\frac d2-\frac dq-\frac\alpha p}\|f\|_{L^2(\mathbb{T}^d)}.
\end{equation}

On the other hand, we consider the eigenfunction $f(x)=\sum_{k\in\mathbb{Z}^d:|k|= N}e^{2\pi ik\cdot x}$. For $R\in\mathbb{N}$, let $$r_d(R)=\# \{k\in\mathbb{Z}^d:|k|^2= R\}.$$ It is well-known that $r_d(R)\ls_\eps R^{\frac{d-2}2+\eps}$ for all $\eps>0$. The bound is essentially sharp, since
\begin{equation}\label{latlow}
	\max_{n\le R}r_d(n)\ge \frac1R\sum_{n\le R}r_d(n)=\frac1R\#\{k\in\mathbb{Z}^d:|k|^2\le R\}\approx R^{\frac{d-2}2}.
\end{equation}

 We have $\|f\|_{L^2(\mathbb{T}^d)}\approx \sqrt{r_d(N^2)}$ and
\[\|e^{it\Delta^{\alpha/2}}f\|_{L^p_{t}L^q_x(\mathbb{T}^{d+1})}\ge	\|f\|_{L^q_x(\{|x|<N^{-1}\})}\approx r_d(N^2) N^{-\frac dq}.\]
This gives a lower bound \begin{equation}\label{lowb2}
	\|e^{it\Delta^{\alpha/2}}f\|_{L^p_{t}L^q_x(\mathbb{T}^{d+1})}\gs \sqrt{r_d(N^2)}N^{-\frac dq}\|f\|_{L^2(\mathbb{T}^d)}.
\end{equation}
Together with \eqref{latlow}, this implies that the lower bound $N^{\frac{d-2}2-\frac dq}$ can be achieved by infinitely many $N$.

For the diagonal cases $p=q$, we restrict to $\alpha >1$. The exponent $\sigma_1(\alpha)$ in Theorem \ref{fracstr} is essentially sharp for the following cases: 
\begin{itemize}
	\item $1<\alpha< 2$ for $q\ge \frac{2(d+2)}d$,
	\item $\alpha=2$ for all $q\ge 2$,
	\item  $\alpha >2$ for $2\le q\le \frac{2(d+2)}d$.
\end{itemize}
For the case $q < \frac{2(d+2)}d$ when $1<\alpha< 2$, the potential for further improvements is tied to weighted $l^2$ decoupling inequalities. Related weighted decoupling inequalities in the weaker $l^p$ setting were recently obtained by Gan-Wu  \cite{gan2025weighted}.

For the sharp Schr\"odinger admissible pair $(p,q)$, we have
\[\frac d2-\frac dq-\frac\alpha p=(2-\alpha)\frac d2(\frac12-\frac1q)=\frac{2-\alpha}p,\]
so the exponent $\sigma_2(\alpha)$ in Theorem \ref{fracstr2} is essentially sharp  when
\begin{itemize}
	\item $0<\alpha<1$ for all $q\ge 2$,
	\item $\alpha>1$ for $2\le q\le \frac{2(d+2)}d$.
\end{itemize}
Moreover, note that $\sigma_2(\alpha)<\sigma_0(\alpha)/2$ when $\alpha\ge2$ with $d>4$ or $1<\alpha<2$ with $d>\frac{2\alpha}{\alpha-1}$, so it improves the exponent in Burq-G\'erard-Tzvetkov \cite{MR2058384}  and  Dinh \cite{dinh17} on general manifolds in these cases.

\subsection{Improved Strichartz estimates for systems}

 To prove Corollary \ref{torussub}, we use Theorem \ref{fracstr2} and Minkowski inequality to get
\begin{equation}\label{str-diag2}
	\Big\|\sum_j  \nu_j |e^{it\Delta^{\alpha/2}} f_j|^2\Big\|_{L^{p/2}_{t}L^{q/2}_x(\mathbb T^{d+1})} \lesssim_{\eps} 
	N^{2\sigma_1+\eps} \|\nu\|_{l^1},
\end{equation}
where $\sigma_1$ is given by \eqref{sigma01}.
Then we obtain \eqref{str-diag} by interpolation between \eqref{str-diag2}, the universal bound \eqref{tb}, and the endpoint estimate \eqref{end0}. 

To prove Corollary \ref{torussup},  by interpolation between the Keel-Tao endpoint estimate \eqref{keeltao} and the kink point estimate \eqref{ktb}, at the Keel-Tao endpoint $(p,q)=(2,\frac{2d}{d-2})$ we have
\begin{equation}\label{ktbound}
	\|	\sum_{j}\nu_j|e^{it\Delta^{\alpha/2}} f_j|^2\|_{L_t^{p/2}L_x^{q/2}}
	\ls N^{\frac2p}\|\nu\|_{l^{\beta}},\ \ \  \forall\beta<\frac{d(d-3)}{(d-2)^2}.
\end{equation}
And then Corollary \ref{torussup} follows from interpolation between \eqref{ktbound} and the subcritical regime in Theorem \ref{thm1}. 

\subsection{Proof of Theorem \ref{dec0}}  Theorem \ref{dec0} can be deduced from the following Lemma \ref{lemmadec} and Minkowski inequality.
\begin{lem}\label{lemmadec}
	Let $\alpha>1$ and $2\le p\le \frac{2(d+2)}{d}$. Then we have for all $\eps>0$,
	\begin{equation}\label{decineq}
		\|E_{[-1,1]^d} g\|_{L^p(\omega_{R_{\delta}})}\ls_\eps \delta^{-\eps}\Big(\sum_{\Delta\in \text{Part}_{\delta^{1/2}}([-1,1]^d)}	\|E_{\Delta} g\|_{L^p(\tilde \omega_{R_{\delta}})}^2\Big)^\frac12
	\end{equation}
	where $R_\delta$ is a rectangular box of size $\delta^{-1}\times ...\times \delta^{-1}\times \delta^{-\max\{1,\alpha/2\}}$.
\end{lem}

\noindent \textbf{Proof.} We  dyadically decompose the cube $[-1,1]^d$ into
\[[-1,1]^d=\{y\in [-1,1]^d:|y|\le \delta^{1/2-\eps}\}\cup \bigcup_{k=1}^K\{y\in [-1,1]^d:2^{k-1}\delta^{1/2-\eps}\le |y|\le 2^k\delta^{1/2-\eps}\}.\]
Here $K\approx \log(\delta^{-1})$. The first part can be easily controlled by Minkowski and H\"older inequalities. It suffices to prove that for any $\delta^{1/2-\eps}\le a\le 1/2$ and the annulus 
\[A_a=\{y\in [-1,1]^d:a\le |y|\le 2a\}\]
we have 
\begin{equation}\label{decineq2}
	\|E_{A_a} g\|_{L^p(\omega_{R_{\delta}})}\ls_\eps \delta^{-\eps}\Big(\sum_{\Delta\in \text{Part}_{\delta^{1/2}}(A_a)}	\|E_{\Delta} g\|_{L^p(\tilde \omega_{R_{\delta}})}^2\Big)^\frac12.
\end{equation}
We claim that for any  $\delta^{1/2-\eps}\le a\le 1/2$, 
\begin{equation}\label{decineq3}
	\|E_{A_a} g\|_{L^p(\omega_{R_{a,\delta}})}\ls_\eps \delta^{-\eps}\Big(\sum_{\Delta\in \text{Part}_{\delta^{1/2}}(A_a)}	\|E_{\Delta} g\|_{L^p( \omega_{R_{a,\delta}})}^2\Big)^\frac12.
\end{equation}
where $R_{a,\delta}$ is a rectangular box of size $\delta^{-1}\times ...\times \delta^{-1}\times a^{2-\alpha}\delta^{-1}$.  Since $\delta^{-\max\{1,\alpha/2\}}\gs a^{2-\alpha}\delta^{-1}$, we can split $R_\delta$ into copies of $R_{a,\delta}$. Note that
\begin{equation}\label{weight}
	\omega_{R_\delta}(x)\ls \sum_{j}\omega_{R_{a,\delta}(j)}(x)\omega_{R_\delta}(j)\ls \tilde\omega_{R_\delta}(x)
\end{equation}
where each $R_{a,\delta}(j)$ is a copy of $R_{a,\delta}$ centered at $j\in \delta^{-1}\mathbb{Z}^d\times a^{2-\alpha}\delta^{-1}\mathbb{Z}$, and the implicit constants only depend on $d$. Then \eqref{decineq3} implies \eqref{decineq2} by 
\eqref{weight} and Minkowski inequality. It is worth to mention that $\sum_{j}\omega_{R_{a,\delta}(j)}(x)\omega_{R_\delta}(j)\approx \omega_{R_\delta}(x)$ cannot hold, so we should use two slightly different rectangular  weight functions in \eqref{weight}, which is different from  the forthcoming \eqref{weight2} for the cube weight functions. 

Now we prove \eqref{decineq3}. Let $\phi(y)=|y|^\alpha$. For $y_0\in A_a$, the Taylor expansion
\[	a^{2-\alpha}\phi(y_0+z)=a^{2-\alpha}(\phi(y_0)+\phi'(y_0)z+\frac12z^T\phi''(y_0)z)+O(a^{-1}|z|^3)\]
where 
\[\phi''(y_0)=\alpha|y_0|^{\alpha-2}(I_d+(\alpha-2)\frac{y_0y_0^T}{|y_0|^2}).\]
Thus, $a^{2-\alpha}|z^T\phi''(y_0)z|\approx |z|^2$ whenever $\alpha>1$.

When $|z|\le \delta^{\frac12-\sigma}$ with $\sigma=\eps/3$, the error term $O(a^{-1}|z|^3)=O(\delta)$ since $a\ge \delta^{1/2-\eps}$. So when $|y-y_0|\le \delta^{\frac12-\sigma}$, the surface $S_a(y)=(y,a^{2-\alpha}|y|^\alpha)$ is in the $\delta$-neighborhood of the paraboloid
\begin{equation}\label{para}
	\tilde S_a(y)=(y, a^{2-\alpha}(\phi(y_0)+\phi'(y_0)(y-y_0)+\frac12(y-y_0)^T\phi''(y_0)(y-y_0))).
\end{equation}

Now we rescale the last variable. Let 
\[E_{\Delta,S_a} g(x)=\int_\Delta g(y)e(x_1y_1+...+x_dy_d+a^{2-\alpha}x_{d+1}|y|^\alpha )dy\]
Then
\[E_{\Delta}g(x)=E_{\Delta,S_a}g(x_1,...,x_d,a^{2-\alpha}x_{d+1})\]
and 
\[\|E_\Delta g\|_{L^p(R_{a.\delta})}=a^{\frac{2-\alpha}p}\|E_{\Delta,S_a} g\|_{L^p(Q_\delta)}\]
where $Q_\delta$ is a cube of side length $\delta^{-1}$. By the scaling, to prove \eqref{decineq3}, it suffices to show
\begin{equation}\label{decineq4}
	\|E_{A_a,S_a} g\|_{L^p(\omega_{Q_\delta})}\ls_\eps \delta^{-\eps}\Big(\sum_{\Delta\in \text{Part}_{\delta^{1/2}}(A_a)}	\|E_{\Delta,S_a} g\|_{L^p( \omega_{Q_\delta})}^2\Big)^\frac12.
\end{equation}
It suffices to show that the smallest constant $K_p(\delta)$ that makes the following inequality holds satisfies $K_p(\delta)\ls_\eps \delta^{-\eps}$
\begin{equation}\label{decineq5}
	\|E_{A_a,S_a} g\|_{L^p(\omega_{Q_\delta})}\le K_p(\delta)\Big(\sum_{\Delta\in \text{Part}_{\delta^{1/2}}(A_a)}\|E_{\Delta,S_a} g\|_{L^p(\omega_{Q_\delta})}^2\Big)^\frac12
\end{equation}
Note that
\begin{equation}\label{weight2}
	\sum_{j}\omega_{Q_{\delta^{1-2\sigma}}(j)}(x)\omega_{Q_\delta}(j)\approx \omega_{Q_\delta}(x)
\end{equation}
where  each $Q_{\delta^{1-2\sigma}}(j)$ is a copy of $Q_{\delta^{1-2\sigma}}$ centered at $j\in \delta^{-1+2\sigma}\mathbb{Z}^{d+1}$, and the implicit constants only depend on $d$. 
By \eqref{decineq5} and Minkowski inequality, we get
\[		\|E_{A_a,S_a} g\|_{L^p(\omega_{Q_\delta})}\le C K_p(\delta^{1-2\sigma})\Big(\sum_{\tau\in \text{Part}_{\delta^{\frac12-\sigma}}(A_a)}\|E_{\tau,S_a} g\|_{L^p(\omega_{Q_\delta})}^2\Big)^\frac12\]Applying Bourgain-Demeter's decoupling inequality \cite[Theorem 1.1]{MR3374964} in the weighted version  \cite[Prop. 9.15]{MR3971577} to the paraboloid \eqref{para}, we get
\[\|E_{\tau,S_a} g\|_{L^p(\omega_{Q_\delta})}\le D_p(\delta)\Big(\sum_{\Delta\in \text{Part}_{\delta^{1/2}}(A_a),\Delta\subset \tau}\|E_{\Delta,S_a} g\|_{L^p(\omega_{Q_\delta})}^2\Big)^\frac12\]
where the decoupling constant $D_p(\delta)\le C_\eps \delta^{-\frac12\eps^2}$ for all $\eps>0$. Thus,
\[K_p(\delta)\le CD_p(\delta)K_p(\delta^{1-2\sigma}).\]Recall $\sigma=\eps/3$. We iterate to get
\begin{align*}
	K_p(\delta)&\le C^k D_p(\delta)D_p(\delta^{1-2\sigma})...D_p(\delta^{(1-2\sigma)^{k-1}})K_p(\delta^{(1-2\sigma)^k})\\
	&\le C^kC_\eps^k\delta^{-\frac12\eps^2(1+(1-2\sigma)^2+...+(1-2\sigma)^{k-1})}K_p(\delta^{(1-2\sigma)^k})\\
	&=C^kC_\eps^k\delta^{-\frac34\eps (1-(1-2\sigma)^k)}K_p(\delta^{(1-2\sigma)^k}).
\end{align*}
Recall that $\delta^{1/2-\eps}\le a\le 1/2$. We choose $k$ such that $\delta^{(1-2\sigma)^k}\approx a^2\le 1/4$, then $K_p(\delta^{(1-2\sigma)^k})\approx 1$ and $k\ls \log\log (\delta^{-1})$. Thus
\[K_p(\delta)\ls C^kC_\eps^k\delta^{-\frac34\eps (1-(1-2\sigma)^k)}\ls_\eps \delta^{-\eps}.\]

\section{Necessary conditions}\label{mr}
Let $d\ge 1,\,\alpha\in (0,\infty)\setminus\{1\},\ N\ge10$. Let $P=\Delta^{\alpha/2}$ or $\sqrt{m^2+\Delta}$ with $m\ge0$.  We make  several crucial observations that are used to show the optimality of the range of $\beta$ in our main theorems. We shall use the Weyl law and the zonal spherical harmonics to construct examples. As before, we fix $(e_k)_k$ to be an orthonormal eigenbasis in $L^2(M)$ associated with the eigenvalues $(\la_k)_k$ of $\sqrt{\Delta}$. Here $0=\la_0<\la_1\le \la_2\le ...$ are arranged in increasing order and we account for multiplicity. In this section, we also discuss the dependence on the interval $I$ so we do not assume $|I|\approx 1$.

$ $

\noindent\textbf{Observation 1.} For all $p,q\ge2$, we have the universal bound  
\begin{equation}\label{tb}
	\Big\|\sum_j  \nu_j |e^{itP} f_j|^2\Big\|_{L_t^{p/2}L_x^{q/2}(I\times M)} \lesssim |I|^{\frac2p}N^{d}
	\|\nu\|_{l^\infty}\end{equation}
for  all orthonormal systems $(f_j)_j$ in $L^2(M)$ with $supp\, \hat f_j\subset \{k:\la_k\le N\}$ and all sequences $\nu=(\nu_j)_j\in l^\beta$. Indeed, for each $t$ the system $\{e^{itP} f_j\}_j$ is also orthonormal  in $L^2(M)$ with Fourier coefficients supported in $\{k:\la_k\le N\}$, then we have 
\begin{align*}
	\sum_j|e^{itP} f_j(x)|^2\le \sum_{k:\la_k\le N}|e_k(x)|^2\ls N^d,\ \ \ \forall t\in I,\ \forall x\in M
\end{align*}
by the pointwise Weyl law.
Moreover, the exponent of $N$ in \eqref{tb} cannot be replaced by any number less than $d$, since if we fix  $\nu_j=1$ for $j\in\{k:\la_k\le N\}$ then	by the Weyl law 
\begin{align}\label{weyl}
	\Big\|\sum_j \nu_j|e^{itP} f_j|^2\Big\|_{L_x^1(M)}\approx \sum_{k:\la_k\le N} 1\approx N^d.
\end{align}

$ $

\noindent\textbf{Observation 2.} The condition 
\begin{equation}\label{nec} \frac1\beta\ge  1-\frac\sigma d \end{equation} is necessary for the estimate \begin{equation}\label{ext-sharp}
	\Big\|\sum_j  \nu_j |e^{itP}f_j|^2\Big\|_{L^{p/2}_tL^{q/2}_x(I\times M)} \lesssim |I|^{\frac2p}N^{\sigma} \|\nu\|_{l^\beta}\end{equation}
to hold for  all orthonormal systems $(f_j)_j$ in $L^2(M)$ with $supp\, \hat f_j\subset \{k:\la_k\le N\}$, and all sequences $\nu=(\nu_j)_j\in l^\beta$. Indeed, 
if we fix  $\nu_j=1$ with $j\in\{k:\la_k\le N\}$ then by \eqref{weyl} we see that \eqref{ext-sharp} implies \eqref{nec}.

$ $

\noindent\textbf{Observation 3.}
Let  $M=S^d$ be the standard sphere, and we fix a point $x_0\in M$. Recall that each eigenvalue of $\sqrt{\Delta}$ on the sphere has the form $\mu_j=\sqrt{j(j+d-1)}$ with $j\in \mathbb{N}$, and for each $\mu_j\approx N$ we can find an $L^2$-normalized zonal function $Z_j$ with $Z_j(x)\approx N^{\frac{d-1}2}$ when the distance $d_g(x,x_0)\ll N^{-1}$. See e.g. Sogge \cite{SoggeHangzhou}.  These $Z_j$ form an  orthonormal system of size $\approx N$ since they are associated with distinct eigenvalues. Since $|e^{itP}Z_j|=|Z_j|$,  we have
\[	\Big\|\sum_j   |e^{itP}Z_j|^2\Big\|_{L^{p/2}_tL^{q/2}_x(I\times M)} \gs N^{d} (N^{-d})^{\frac2q}|I|^{\frac2p}.\]
This gives another necessary condition for \eqref{ext-sharp} on the sphere 
\begin{equation}\label{nec2}
	\frac1\beta\ge d-\frac{2d}q-\sigma.
\end{equation}
If the sharp Schr\"odinger admissible condition $\frac1p=\frac d2(\frac12-\frac1q)$ holds, then \eqref{nec2} is equivalent to 
\begin{equation}\label{nec22}
	\frac1\beta\ge \frac4p-\sigma.
\end{equation}

\noindent\textbf{Observation 4.} Let $M$ be a compact manifold of dimension $d$. Fix any point $x_0\in M$. By the pointwise Weyl law, we have for all large $j$
\[\sum_{k:\la_k\le j}|e_k(x_0)|^2=C_d j^d+O(j^{d-1}).\]
 Then  there exists a constant $c>0$ such that for all large $j$
\[\sum_{k:\la_k\in(j-c,j]}|e_k(x_0)|^2\approx  j^{d-1}.\]
So we can construct an orthonormal system $(f_j)_j$ of size $\approx N$, where $j\approx N$ and
\[f_{j}(x)=c_j\sum_{k:\la_k\in(j-c,j]}e_k(x_0)\overline{e_k(x)}.\]
Here the normalized constant $c_j=(\sum_{k:\la_k\in(j-c,j]}|e_k(x_0)|^2)^{1/2}\approx j^{\frac{d-1}2}.$
So $f_j(x)\approx N^{\frac{d-1}2}$ when $d_g(x,x_0)\ll N^{-1}$. Indeed, this follows from the mean value theorem and Bernstein inequality
\[\|\nabla f_j\|_{L^\infty(M)}\ls N\|f_j\|_{L^\infty(M)}.\] See e.g. Imekraz- Ouhabaz \cite{IO22}. Then
\[e^{it\Delta^{\alpha/2}}f_j=c_je^{itj^\alpha}\sum_{k:\la_k\in(j-c,j]} e^{it(\la_k^\alpha-j^\alpha)}e_k(x_0)\overline{e_k(x)}.\]
Note that 
\[|(j-c)^\alpha-j^\alpha|\ls j^{\alpha-1}\approx N^{\alpha-1}.\]
So  we still have $|e^{it\Delta^{\alpha/2}}f_j(x)|\approx N^{\frac{d-1}2}$ when $d_g(x,x_0)\ll N^{-1}$ and $|t|\ll N^{1-\alpha}$. Thus, for $|I|\gs N^{1-\alpha}$ we have
\[	\Big\|\sum_j   |e^{it\Delta^{\alpha/2}}f_j|^2\Big\|_{L^{p/2}_tL^{q/2}_x(I\times M)} \gs N^{d} (N^{-d})^{\frac2q} N^{\frac2p(1-\alpha)}.\]
When $0<\alpha\le1$ and $|I|\approx N^{1-\alpha}$, this gives a necessary condition for \eqref{ext-sharp} on any compact manifold 
\begin{equation}\label{nec3}
	\frac1\beta\ge d-\frac{2d}q-\sigma.
\end{equation}
If the sharp Schr\"odinger admissible condition $\frac1p=\frac d2(\frac12-\frac1q)$ holds, then \eqref{nec3} is equivalent to 
\begin{equation}\label{nec221}
	\frac1\beta\ge \frac{4}p-\sigma.
\end{equation}
If  the sharp wave admissible condition $\frac1p=\frac {d-1}2(\frac12-\frac1q)$ holds, then \eqref{nec3} is equivalent to 
\begin{equation}\label{nec222}
	\frac1\beta\ge \frac{d}{d-1}\frac4p-\sigma.
\end{equation}

\section{Applications}
As in the works by Frank–Sabin \cite{MR3730931}, Lewin–Sabin \cite{MR3304272,MR3270166}, Nakamura \cite{MR4068269}, Bez-Lee-Nakamura  \cite{BLN21},  we can exploit the Strichartz estimates to prove the well-posedness of  the infinite systems of dispersive equations  with Hartree-type nonlinearity on compact manifolds
\begin{equation}\label{hart}
	\begin{cases}
		i\partial_t u_j=P u_j+(W\rho)u_j,\ \ j\in\mathbb{N}\\
		u_j(0,\cdot)=f_j
	\end{cases}
\end{equation}
where $\rho=\sum_{j=1}^\infty |u_j|^2$ and $W\rho$ is a real-valued function on $M$. We focus on $P=\Delta^{\alpha/2}$ or $\sqrt{m^2+\Delta}$ with $m\ge0$. In the flat case, it is standard to take the convolution operator $W \rho=w*\rho$, where $w$ is the interaction potential function on $M$. 

\subsection{Conditions on the systems}
Let $s\ge0$. First, we need the Strichartz estimates
\begin{equation}\label{ass1}
	\Big\|\sum_j  \nu_j |e^{itP} f_j|^2\Big\|_{L_t^{p/2}L_x^{q/2}([0,1]\times M)} \ls
	\|\nu\|_{l^\beta}\end{equation}
 for all orthonormal systems $(f_j)_j$ in $H^{s}(M)$ and all sequences $\nu=(\nu_j)_j\in l^\beta$. As we have seen in the introduction, the estimates \eqref{ass1} can be deduced from the frequency localized estimates 
\eqref{eq0} as in Theorem \ref{thm1}, \ref{thm2},  \ref{wkgthm}, \ref{wkgthm2} and Corollary \ref{torussub}, \ref{torussup}. Moreover, \eqref{ass1} implies for any $T>0$,
\[	\Big\|\sum_j  \nu_j |e^{itP} f_j|^2\Big\|_{L_t^{p/2}L_x^{q/2}([0,T]\times M)} \ls
T^{\frac2p}\|\nu\|_{l^\beta}.\]

Let $\D=\sqrt{1+\Delta}$. We also need the control of the Hartree-type nonlinearity
\begin{equation}\label{ass2}
	\|\D^{\pm s}(W\rho)\D^{\mp s}\|_{\mathfrak S^\infty}\le C_{s,q,W}\|\rho\|_{L^{q/2}(M)}
\end{equation}
or equivalently
\[	\|(W\rho)f\|_{H^{ r}(M)}\le C_{s,q,W}\|\rho\|_{L^{q/2}(M)}\|f\|_{H^{r}(M)},\ \forall f\in H^r(M),\ \ r=\pm s.\]
In the flat case with $W\rho=w*\rho$, the condition \eqref{ass2}  holds true with $C_{s,q,W}=C_{s,\delta}\|w\|_{B_{(q/2)',\infty}^{s+\delta}}$ for all $\delta>0$ (see \cite{MR4068269, BLN21}). Indeed, it follows from the inequality for the Besov norm
 \[	\|gf\|_{H^r}\le C_{r,\delta}\|g\|_{B_{\infty,\infty}^{|r|+\delta}}\|f\|_{H^r},\ \ \forall r\in \mathbb{R},\ \forall \delta>0,\]
 and by H\"older inequality
\[	\|w*\rho\|_{B_{\infty,\infty}^{s+\delta}}\le \|w\|_{B_{(q/2)',\infty}^{s+\delta}} \|\rho\|_{L^{q/2}}\ \ .\]
See Triebel \cite[p. 29 \& p. 205]{MR1163193} and Seeger-Sogge \cite[Theorem 4.1]{ss89} for characterizations of Besov space on compact manifolds. A  typical example in the flat case is $w(x)=|x|^{-a}$ with $a<d$. See Nakamura \cite{MR4068269}. A natural generalization of  this convolution operator on compact manifolds is the spectral multiplier $\D^{-d+a}$. We calculate the norm
\begin{align*}
	\|\D^{-d+a}\rho\|_{B_{\infty,\infty}^s}&=\sup_{j\ge0}2^{js}\|\varphi_j(\D)\D^{-d+a}\rho\|_{L^\infty}\\
	&\ls \sup_{j\ge0}2^{j(s-d+a)}\|\tilde\varphi_j(\D)\rho\|_{L^\infty}\\
	&\ls \sup_{j\ge0}2^{j(s-d+a+2d/q)}\|\rho\|_{L^{q/2}},
\end{align*}
where $\tilde\varphi_j$ shares essentially the same property as the Littlewood-Paley bump function $\varphi_j$.
So we require that $s-d+a+2d/q\le 0$, which is equivalent to $a\le d-\frac{2d}q-s$. For the sharp Schr\"odinger admissible pairs $(p,q)$, it is equivalent to $a\le \frac4p-s$. For the sharp wave admissible pairs $(p,q)$,  it is equivalent to $a\le \frac{4d}{p(d-1)}-s$.

\subsection{Well-posedness of the systems}
For applications, it is useful to state the condition \eqref{ass1} in an operator-theoretic version. Given a compact self-adjoint operator $\gamma$ on $L^2(M)$, by the spectral theorem we can write 
\[\gamma h=\sum_j \nu_j \langle h,f_j\rangle f_j,\ \ \forall h\in L^2(M).\]We formally denote the diagonal of the integral kernel of $\gamma$ by
\[\rho_\gamma(x)=\sum_j\nu_j|f_j(x)|^2.\]
By the assumption \eqref{ass1}, $\rho_{\D^{-s}\gamma(t)\D^{-s}}$ is well-defined in $L^{p/2}_tL^{q/2}_x$, where $\gamma(t)=e^{itP}\gamma_0e^{-itP}$, and satisfies
\[\|\rho_{\D^{-s}\gamma(t)\D^{-s}}\|_{L^{p/2}_tL^{q/2}_x}\le C_*\|\gamma_0\|_{\mathfrak S^\beta}\]
whenever $\gamma_0\in \mathfrak S^\beta$. We define the Sobolev-type Schatten norm by
\[\|\gamma\|_{\mathfrak S^{\beta,s}}=\|\D^s\gamma\D^s\|_{\mathfrak S^{\beta}}.\]
It is standard to transform the infinite system \eqref{hart} into the operator formalism
\begin{equation}\label{opeq}
	\begin{cases}
		i\partial_t \gamma=[P+W\rho_\gamma,\gamma],\ \ (t,x)\in \mathbb{R}\times M,\\
		\gamma(0,\cdot)=\gamma_0.
	\end{cases}
\end{equation}
The following local well-posedness was proved in the abstract form in \cite[Prop. 10]{BLN21}, and it still holds on compact manifolds. The proof is based on the Duhamel principle and the contraction mapping theorem.
\begin{prop}[Local well-posedness]\label{blnprop}
	Suppose \eqref{ass1} and \eqref{ass2} hold. Then for any $\gamma_0\in \mathfrak S^{\beta,s}$, there exist $T=T(\|\gamma_0\|_{\mathfrak S^{\beta,s}},C_{s,q,W})>0$ and a unique solution $\gamma\in C_t^0([0,T];\mathfrak S^{\beta,s})$ to \eqref{opeq} on $[0,T]\times M$ with $\rho_\gamma\in L^{p/2}_tL^{q/2}_x$.
\end{prop}
We also have the following  global well-posedness for small data. The argument is similar, see \cite[Prop. 4.1]{MR4068269}. 
\begin{prop}[Almost global well-posedness]\label{blnprop2}
	Suppose \eqref{ass1} and \eqref{ass2} hold. Then each $T>0$,  there exists $\delta_T=\delta_T(C_{s,q,W})$ such that for any for any $\|\gamma_0\|_{\mathfrak S^{\beta,s}}<\delta_T$, there exists a unique solution $\gamma\in C_t^0([0,T];\mathfrak S^{\beta,s})$ to \eqref{opeq} on $[0,T]\times M$ with $\rho_\gamma\in L^{p/2}_tL^{q/2}_x$.
\end{prop}
As corollaries, we obtain the  well-posedness of \eqref{opeq} for the Hartree-type nonlinearity  $W\rho=\D^{-d+a}\rho$ by Theorem \ref{thm1} and Theorem \ref{wkgthm}.
\begin{cor}
	Let $d\ge1$ and $\alpha>1$. Suppose $(p,q)$ and $\beta$ are as in Theorem \ref{thm1}. Let $s>\frac1p$ and $a\le \frac4p-s$. Let $P=\Delta^{\alpha/2}$ and $W\rho=\D^{-d+a}\rho$.  Then the system \eqref{opeq} has local well-posedness and almost global well-posedness as in Propositions \ref{blnprop} and \ref{blnprop2}.
\end{cor}
\begin{cor}
	Let $d\ge2$. Suppose $(p,q)$ and $\beta$ are as in Theorem \ref{wkgthm}. Let $s>\frac1p\frac{d+1}{d-1}$ and $a\le \frac{4d}{p(d-1)}-s$. Let $P=\sqrt{m^2+\Delta}$ with $m\ge0$ and $W\rho=\D^{-d+a}\rho$.  Then the system \eqref{opeq} has local well-posedness and almost global well-posedness as in Propositions \ref{blnprop} and \ref{blnprop2}.
\end{cor}

\section*{Declarations}
\noindent \textbf{Data availability statement.} Data sharing not applicable to this article as no datasets were generated or analyzed during the current study.

\noindent \textbf{Conflict of interests.} The authors have no relevant financial or non-financial interests to disclose.


	\bibliography{os}
	
	\bibliographystyle{plain}

\end{document}